\documentclass[11pt,reqno,letter]{amsart}

\usepackage{amsmath}
\usepackage{amssymb}
\usepackage[dvipdfmx]{graphicx}
\usepackage{hyperref}
\usepackage[usenames]{color}
\usepackage{caption}
\usepackage{subcaption}
\usepackage{pifont}
\usepackage{bm}
\usepackage{url}
\usepackage{setspace}
\usepackage[margin=1in]{geometry}
\usepackage{subcaption}
\usepackage{xfrac}
\usepackage{appendix}
\usepackage[numbers]{natbib}
\usepackage{algorithm}
\usepackage{algpseudocode}

\usepackage[normalem]{ulem}

\usepackage{float}

\theoremstyle{plain}
\newtheorem{theorem}{Theorem}[section]

\newtheorem{corollary}[theorem]{Corollary}
\newtheorem{proposition}[theorem]{Proposition}

\theoremstyle{definition}

\newtheorem{remark}[theorem]{Remark}

\numberwithin{equation}{section}
\numberwithin{theorem}{section}
\numberwithin{table}{section}
\numberwithin{figure}{section}

\def\@fnsymbol#1{\ensuremath{\ifcase#1\or *\or \dagger\or \ddagger\or
   \mathsection\or \mathparagraph\or \|\or **\or \dagger\dagger
   \or \ddagger\ddagger \else\@ctrerr\fi}}

\usepackage{hyperref}%

\title[EM algorithms for optimization problems]{EM algorithms for optimization problems \\
with polynomial objectives}

\author[Asai]{Kensuke Asai\textsuperscript{$1*\dagger$}}
\thanks{
{\it Disclaimer}: The first author of this paper hereby states that his affiliated institution has had no direct involvement in the research and bears no responsibility for the findings or conclusions of this study.}
\author[Gotoh]{Jun-ya Gotoh\textsuperscript{$2\dagger$}}

\dedicatory{
\textsuperscript{$1*$}
{Coincheck, Inc.,} Tokyo, Japan. \\
\textsuperscript{$2$}Department of Data Science for Business Innovation, Chuo University, Tokyo, Japan. \\
\, \\
\textsuperscript{$*$}
Corresponding author(s). E-mail(s): kensuke.asai.quants@gmail.com \\
Contributing authors: jgoto@kc.chuo-u.ac.jp \\
\textsuperscript{$\dagger$}
These authors contributed equally to this work.
}

\begin{document}

\begin{abstract}
The EM (Expectation--Maximization) algorithm is regarded as an MM (Majorization--Minimization) algorithm for maximum likelihood estimation of statistical models. Expanding this view, this paper demonstrates that by choosing an appropriate probability distribution, even nonstatistical optimization problem can be cast as a negative log-likelihood-like minimization problem, which can be approached by an EM (or MM) algorithm. When a polynomial objective is optimized over a simple polyhedral feasible set and an exponential family distribution is employed, the EM algorithm can be reduced to a natural gradient descent of the employed distribution with a constant step size. This is demonstrated through three examples. In this paper, we demonstrate the global convergence of specific cases with some exponential family distributions in a general form. In instances when the feasible set is not sufficiently simple, the use of MM algorithms can nevertheless be adequately described. When the objective is to minimize a convex quadratic function and the constraints are polyhedral, global convergence can also be established based on the existing results for an entropy-like proximal point algorithm.
\end{abstract}

\maketitle
%
%

\section{Introduction}
After being rediscovered many times by predecessors, the EM (Expectation--Maximization) algorithm was coined and popularized by Dempster, Laird, and Rubin \cite{dempster1977maximum} in late 1970s as a framework used in estimating statistical models, in particular, when those depend on latent random variables. 
See, for example, 
\cite{mclachlan2007algorithm} for plenty of examples of successful applications of the framework and extensive references therein. In the context of mathematical optimization, the EM algorithm can be viewed as a special example of MM (Majorization--Minimization) algorithm (e.g., \cite{hunter2004tutorial,lange2016mm}) or a proximal-like algorithm that employs the Kullback--Leibler divergence (e.g., \cite{tseng2004analysis}). Focusing on such an MM aspect, we study solving a general, i.e., not necessarily statistical, optimization problem via an EM algorithm (without observed data) by rewriting it to another optimization of log-likelihood-like function with the help of 
a certain probability distribution and by reading the parameters of the distribution as the decision variables of optimization problems in hand. 

Table \ref{tbl:problems} summarizes examples of such combination of optimization problem and exponential family distribution explored in this article. In particular, for the first three cases (i), (ii), and (iii) in Table \ref{tbl:problems}, we show that the EM algorithm results in a gradient descent algorithm with a natural gradient \cite{amari1996neural} with a constant step size. Even when the constraints are not simple, as (iv) to (vi), the EM scheme yields MM algorithms which can further derive EM gradient algorithms \cite{lange1995gradient}, under which a bit more practical presentations of interior point methods can be derived. 

The contribution of this article is to bridge a few areas whose connections have been less recognized because of different terminology. 
On the other hand, it is not our intention to present a cutting-edge efficient algorithm for optimization problems, but a unified framework over which new developments are expected to be made in future researches.  

\begin{table}[h]
\centering
\caption{Problem types and used probability distributions for EM algorithms}
\label{tbl:problems}
\begin{tabular}{ccccc}
\hline
      & objective & constraints & probability & natural grad. descent\\
      & function &                 & distribution & w/ const. step size\\
\hline
(i)   & convex quadratic & --- & normal (const.cov.) & yes \\
(ii)  & polynomial & rectangle & binomial & yes \\
(iii) & polynomial & simplex & multinomial & yes \\
(iv) & polynomial & polytope & binomial & no \\
(v)  & linear & polyhedron & Poisson & no \\
(vi) & 
quadratic & polyhedron & Poisson$+$normal (const.cov.)& no \\
\hline
\end{tabular}
\end{table}

\paragraph{\it Relations to existing works}

Researchers have employed the optimization of parameters of probability mass functions belonging to the exponential family in contexts other than purely statistical estimation, such as maximum likelihood estimation. 
For example, \cite{shirakawa2018dynamic} presented a methodology for automatically estimating the structure of a neural network by minimizing the expected loss under a multivariate Bernoulli distribution, where the optimization variables are the parameters of the distribution, using a stochastic gradient descent method based on the natural gradient. 
In contrast to those studies, this paper does not specify an application. Rather, it systematically demonstrates that by selecting an appropriate exponential distribution family, it is possible to construct an EM (MM) algorithm and a natural gradient descent method, even for a general class of optimization problems that do not contain any statistical elements.

In terms of focusing on the MM aspect of the EM algorithm and extending it to general optimization problems, the paper by Becker et al. \cite{becker1997algorithms} is closer to ours, even though their argument starts from the statistical EM algorithm. 
To make a similar claim, they present the EM algorithms in the absence of missing data, which is typically assumed when working with the EM algorithms. 
In contrast to \cite{becker1997algorithms} and other statistical EM algorithms, this paper focuses on the EM algorithms with no observed data, i.e., all data is missing, to extend to non-statistical optimization problems. This contrast stems from the difference in stand points that \cite{becker1997algorithms} is aimed at statisticians, while this paper is aimed at a broader readership including optimization researchers.
In fact, the present paper puts forth a novel perspective in which an optimization problem devoid of statistical elements can be artificially linked to log-likelihood-like optimization, and through the EM algorithm for that, we derive an MM algorithm for the original optimization problem. 

Some existing works have proposed interior point algorithms for linearly constrained  convex program, employing a logarithmic barrier term that appears in the context of EM algorithm 
(e.g., \cite{eggermont1990multiplicative,lange1994adaptive,iusem1995convergence,auslender1995interior,teboulle1997convergence}). 
We show that an MM algorithm for quadratic program (vi) (and linear program as a special case) can be constructed in the combination with a hybrid of the Poisson and normal distributions and is a proximal-like algorithm with a hybrid proximal term; 
that for linear program (v) can also be constructed combined only with the Poisson and is a special case of \cite{teboulle1997convergence}, which considers a proximal-like algorithm with $\phi$-divergence-based proximal term. The global convergence of our MM algorithm for convex QP and LP can be proven based on \cite{teboulle1997convergence}. 

 For the nonnegativity-constrained least squares, a gradient projection algorithm based on a quadratic surrogate function has been known (see \cite{sun2016majorization}), and the surrogate function can be improved under additional assumptions on the design matrix \cite{de1987convergence,lee2000algorithms}. 
In contrast, we can present a couple of different constructions of surrogate functions for the same problem by employing different probability distributions, each leading to an interior point algorithm. 

Lastly, the utilization of probability models for optimization may evoke associations with model-based random search methods (see, e.g., \cite{hu2012survey,bartz2017model} for a survey). For example, exponential family distributions have been used in the Cross Entropy (CE) method \cite{rubinstein1999cross,de2005tutorial} for combinatorial and nonconvex continuous optimization problems. 
While those methods correspond the decision variables to the random variables and conduct a stochastic search, our approach is purely deterministic procedure by corresponding the decision variables to the parameters. 

\paragraph{\it Structure of the paper}
The next section describes basic strategy and required components of the proposed scheme. Section \ref{sec:expEM} shows properties of the scheme for the case where the involved distribution is a member of the exponential family. 
Section \ref{sec:examples} presents three examples (i) to (iii) in Table \ref{tbl:problems}. Section \ref{sec:further_examples} is devoted to further examples (iv) and (vi), which are based on shifted exponential family distributions. Section \ref{sec:conclusion} concludes the paper. The majority of the proofs of propositions are located in the appendix, with the intention of enhancing the readability of the manuscript. 
The example (v) is also described in the appendix. 
\paragraph{\it Notation} 
We use small bold letters for indicating vectors and capital bold letters for indicating matrices. In particular, $\bm{0},\bm{1}$ denote the all-zero or all-one vectors, respectively. Small italic letters are used to stand for components of vectors and matrices, e.g., $\bm{\theta}=(\theta_1,...,\theta_p)^\top$. 
$\mathbb{R}^p_+$ and $\mathbb{R}^p_{++}$ denote the spaces of nonnegative and positive vectors, respectively.
We denote by $\mathbb{S}^p$ the space of $p\times p$ real symmetric matrices, by  
$\mathbb{S}^p_+,\mathbb{S}^p_{++}$ the spaces of positive semidefinite and positive definite matrices in $\mathbb{S}^p$, respectively. For $\bm{A},\bm{B}\in\mathbb{S}^p$, $\bm{A}\succ\bm{B}$ stands for $\bm{A}-\bm{B}\in\mathbb{S}^p_{++}$. 
By $\|\cdot\|$ we denote the Euclidean norm. For $\bm{Q}\in\mathbb{S}^p_{++}$, we denote the ellipsoidal norm by $\|\cdot\|_{\bm{Q}}$, i.e., $\|\bm{z}\|_{\bm{Q}}:=\sqrt{\bm{z}^\top\bm{Q}\bm{z}}$. 
We use the abbreviation $[p]:=\{1,...,p\}$.


\section{EM scheme for general optimization problem}\label{sec:framework}
We consider minimizing a function $f$ over a closed convex set $\Theta\subset\mathbb{R}^p$: 
\begin{align}
\underset{\bm{\theta}\in \Theta}{\mbox{minimize}}~ & f(\bm{\theta}) \label{eq:P}
\end{align}
via 
\emph{MM algorithm}, which alternately minimizes a surrogate function $u(\bm{\theta}|\bar{\bm{\theta}})$ of $f(\bm{\theta})$ such that 
\begin{align}
u(\bm{\theta}\vert\bar{\bm{\theta}}) &\geq f(\bm{\theta}),\quad\mbox{for any }\bm{\theta},\bar{\bm{\theta}}\in\mathrm{int}\Theta
,\label{eq:u.geq.f}\\
u(\bar{\bm{\theta}}\vert\bar{\bm{\theta}}) &= f(\bar{\bm{\theta}}),\quad \mbox{for any }
\bar{\bm{\theta}}\in\mathrm{int}\Theta.\label{eq:u=f}
\end{align}
To define surrogate functions motivated by the EM algorithm, we first introduce random variables and probability distribution that satisfy a certain condition. 

Let $\bm{X}$ be a vector of random variables that has a probability density $p_{\bm{x}}$ such that $\mathbb{P}\{\bm{X}\leq\bm{x}\vert\bm{\theta}\}=
\int_{\bm{\xi}\leq\bm{x}}p_{\bm{\xi}}(\bm{\theta})\mathrm{d}\bm{\xi}$, a probability mass $p_{\bm{x}}(\bm{\theta})=\mathbb{P}\{\bm{X}=\bm{x}\vert\bm{\theta}\}$, or their hybrid, 
where $\bm{\theta}\in\Theta\subset \mathbb{R}^p$ is the vector of parameters of the probability distribution and $\Xi$ is the support of $\bm{X}$. Assume that for \eqref{eq:P} and $p_{\bm{x}}$ satisfying $\mathrm{int}\Theta=\{\bm{\theta}\in\mathbb{R}^p\vert p_{\bm{x}}(\bm{\theta})> 0\}$, there is a positive-valued function $G:\Xi\to(0,\infty)$ satisfying 
\begin{align}
f(\bm{\theta})=-\ln\mathbb{E}_{p_{\bm{x}}(\bm{\theta})}\big[G(\bm{X})\big]\quad\mbox{for all }\bm{\theta}\in
\mathrm{int}\Theta,\label{eq:G-condition}
\end{align}
where 
$\mathbb{E}_{p_{\bm{x}}(\bm{\theta})}\big[\cdot\big]$ denotes the expectation under the distribution $p_{\bm{x}}(\bm{\theta})$. (For the expectation computation in this article, we only describe the case where $p_{\bm{x}}$ is a probability density because it is straightforward to get results for the other cases, for example, by replacing the integral with a discrete summation when $p_{\bm{x}}$ is a probability mass.) 
While we consider minimizing the negative log-likelihood-like function \eqref{eq:G-condition} for \eqref{eq:P}, it corresponds to the EM algorithm where no data point is observed and only the distribution of latent variables is assumed, which is in contrast to the ordinary context or \cite{becker1997algorithms}. 
As will be shown below, however, this set-up is enough to link the EM algorithm to \eqref{eq:P}. Some examples of the tuple $(f,\Theta,p_{\bm{x}},G)$ will be presented in Sections \ref{sec:examples} and \ref{sec:further_examples}. 

For a tuple $(f,\Theta,p_{\bm{x}},G)$ satisfying \eqref{eq:G-condition}, we see that
\begin{align} \label{eq:def_u}
u(\bm{\theta}\vert\bar{\bm{\theta}}) &:= - \int_{\bm{x}\in\Xi}q_{\bm{x}}(\bar{\bm{\theta}})\ln\frac{G(\bm{x}) p_{\bm{x}}(\bm{\theta})}{q_{\bm{x}}(\bar{\bm{\theta}})} \mathrm{d}\bm{x}
\end{align}
satisfies \eqref{eq:u.geq.f} since, for another density $q_{\bm{x}}$ (equivalent to $p_{\bm{x}}$), we have for any $\bm{\theta},\bar{\bm{\theta}}\in\mathrm{int}\Theta$,
\begin{align}
f(\bm{\theta})&
\underset{\eqref{eq:G-condition}}{=} 
-\ln\int_{\bm{x}\in\Xi}q_{\bm{x}}(\bar{\bm{\theta}})\frac{G(\bm{x}) p_{\bm{x}}(\bm{\theta})
}{q_{\bm{x}}(\bar{\bm{\theta}})} \mathrm{d}\bm{x} 
\leq
- \int_{\bm{x}\in\Xi}q_{\bm{x}}(\bar{\bm{\theta}})\ln\frac{G(\bm{x}) p_{\bm{x}}(\bm{\theta})
}{q_{\bm{x}}(\bar{\bm{\theta}})} \mathrm{d}\bm{x}
\textcolor{black}{\underset{\eqref{eq:def_u}}{=}} u(\bm{\theta}\vert\bar{\bm{\theta}}),
\label{eq:how_to_majorize}
\end{align}
where the inequality holds due to the Jensen's inequality. Letting $s(\bm{x}|\bm{\theta}):=G(\bm{x}) p_{\bm{x}}(\bm{\theta}),$ we see that the nonnegative gap, denoted $h(\bm{\theta}|\bar{\bm{\theta}})$, between $u(\bm{\theta}\vert\bar{\bm{\theta}})$ and $f(\bm{\theta})$ equals the Kullback--Leibler (KL) divergence $\mathrm{KL}\big(\bm{q}\big\|\bm{p}\big):=\int_{\bm{x}\in\Xi}q_{\bm{x}}\ln\frac{q_{\bm{x}}}{p_{\bm{x}}} \mathrm{d}\bm{x}$ of the distribution 
$q_{\bm{x}}=
q_{\bm{x}}(\bar{\bm{\theta}})
,\bm{x}\in\Xi$, relative to $p_{\bm{x}}=\frac{s(\bm{x}|\bm{\theta})}{\int_{\bm{\xi}} s(\bm{\xi}|\bm{\theta}) \mathrm{d}\bm{\xi}},\bm{x}\in\Xi$, that is,
\begin{align}
h(\bm{\theta}|\bar{\bm{\theta}}):=u(\bm{\theta}|\bar{\bm{\theta}})-f(\bm{\theta}) &= \int_{\bm{x}\in\Xi} q_{\bm{x}}(\bar{\bm{\theta}})\ln\frac{q_{\bm{x}}(\bar{\bm{\theta}})}{\left(\displaystyle\frac{s(\bm{x}|\bm{\theta})}{\int_{\bm{\xi}\in\Xi} s(\bm{\xi}|\bm{\theta}) \mathrm{d}\bm{\xi}}\right)} \mathrm{d}\bm{x}\geq 0.
 \label{eq:KL-nonnegativity}
\end{align}
Therefore, if we set $\bar{\bm{\theta}}$ so that $q_{\bm{x}}(\bar{\bm{\theta}})=\frac{s(\bm{x}|\bar{\bm{\theta}})}{\int_{\bm{\xi}\in\Xi} s(\bm{\xi}|\bar{\bm{\theta}}) \mathrm{d}{\bm{\xi}}}$ is fulfilled, the function $u$ satisfies
\eqref{eq:u=f}. 

For an incumbent solution $\bm{\theta}=\bm{\theta}^{(t-1)}$ at the $t$-th iteration, the EM algorithm updates the probability distribution by the formula: 
\begin{align}
q_{\bm{x}}(\bm{\theta}^{(t-1)})
&\leftarrow\frac{s(\bm{x}|\bm{\theta}^{(t-1)})}{\int_{\bm{\xi}\in\Xi} s(\bm{\xi}|\bm{\theta}^{(t-1)}) \mathrm{d}{\bm{\xi}}}
=\frac{G(\bm{x})}{\mathbb{E}_{p_{\bm{x}}(\bm{\theta}^{(t-1)})}[G(\bm{X})]}p_{\bm{x}}(\bm{\theta}^{(t-1)}), 
\qquad\mbox{for }\bm{x}\in\Xi,
\label{eq:q_x}
\end{align}
which minimizes the KL divergence relative to $\big(\frac{s(\bm{x}|\bm{\theta}^{(t-1)})}{\int_{\bm{\xi}} s(\bm{\xi}|\bm{\theta}^{(t-1)})\mathrm{d}{\bm{\xi}} }\big)_{\bm{x}}$, ensuring 
$
u(\bm{\theta}^{(t-1)}|\bm{\theta}^{(t-1)})=f(\bm{\theta}^{(t-1)}).
$
The surrogate function $u(\bm{\theta}|\bm{\theta}^{(t-1)})$ is then minimized:
\begin{align}
\bm{\theta}^{(t)}\gets\arg\min_{\bm{\theta}\in\Theta} u(\bm{\theta}|\bm{\theta}^{(t-1)}).
\label{eq:thetabaarisminimizer}
\end{align}
The EM algorithm repeats the two steps as in Algorithm \ref{alg:prototypeEM}.
If the minimizer \eqref{eq:thetabaarisminimizer} is found, 
it is a solution whose objective value is no worse than the incumbent since we have 
\begin{align*}
f(\bm{\theta}^{(t)})
\underset{\eqref{eq:KL-nonnegativity}}{\leq} u(\bm{\theta}^{(t)}|\bm{\theta}^{(t-1)})
\underset{\eqref{eq:thetabaarisminimizer}}{\leq} u(\bm{\theta}^{(t-1)}|\bm{\theta}^{(t-1)})
\underset{\eqref{eq:q_x}
}{=} f(\bm{\theta}^{(t-1)}).
\end{align*}
\begin{algorithm}[tb]
\caption{Prototype of 
EM algorithm for optimization problem \eqref{eq:P}}
\label{alg:prototypeEM}
\begin{algorithmic}[1]
\State Let $(f,\Theta,p_{\bm{x}},G)$ 
be a tuple satisfying \eqref{eq:G-condition}, and $\bm{\theta}^{(0)}\in\mathrm{int}\Theta$, and set $t\gets 1$.
\Repeat
	\State 
$\displaystyle q_{\bm{x}}(\bm{\theta}^{(t-1)})
\gets\frac{s(\bm{x}|\bm{\theta}^{(t-1)})}{\int_{\bm{\xi}\in\Xi} s(\bm{\xi}|\bm{\theta}^{(t-1)}) \mathrm{d}\bm{\xi}}=\frac{G(\bm{x}) p_{\bm{x}}(\bm{\theta}^{(t-1)})}{\int_{\bm{\xi}\in\Xi}G(\bm{\xi}) p_{\bm{\xi}}(\bm{\theta}^{(t-1)}) \mathrm{d}\bm{\xi}},\qquad\mbox{for }\bm{x}\in\Xi
$
	\State
$\displaystyle \bm{\theta}^{(t)}\gets\arg\min_{\bm{\theta}\in\Theta}
u(\bm{\theta}|\bm{\theta}^{(t-1)})=
\arg\min_{\bm{\theta}\in\Theta}\Big\{
-\int_{\bm{x}\in\Xi} q_{\bm{x}}(\bm{\theta}^{(t-1)}) \ln s(\bm{x}|\bm{\theta})  \mathrm{d}\bm{x}\Big\}.$
\Until{a termination condition is fulfilled}
\end{algorithmic}
\end{algorithm}
The convergence of the sequence $\{\bm{\theta}^{(t)}\}_{t\geq 0}$ to a stationary point is known under some conditions (e.g., Theorem 1 of \cite{razaviyayn2013unified}). However, those conditions do not hold true for some applications in particular when stationary point is attained at the boundary of the feasible region $\Theta$. 

\begin{remark}
\label{sec:composite}
The above MM scheme can be extended to the composite objective's case:
\begin{align}
\underset{\bm{\theta}\in\Theta}{\mbox{{minimize}}}~ & \Phi(\bm{\theta}) := f(\bm{\theta}) + \rho(\bm{\theta}),
\label{eq:F}
\end{align}
where $\rho$ is a function on $\Theta$. For any surrogate function $u(\cdot\vert\cdot)$ for $f(\cdot)$, 
$
U(\bm{\theta}\vert\bar{\bm{\theta}}) := u(\bm{\theta}\vert\bar{\bm{\theta}}) + \rho(\bm{\theta})
$
is a surrogate function for $\Phi(\bm{\theta})$ and an MM algorithm for \eqref{eq:F} 
can be constructed 
by simply replacing the M-step (line 4) of Algorithm \ref{alg:prototypeEM} with
$
\bm{\theta}^{(t)}\gets\arg\min\limits_{\bm{\theta}\in\Theta}u(\bm{\theta}|\bm{\theta}^{(t-1)})+ \rho(\bm{\theta})
$.
An example of this extension will be given in Section \ref{sec:l1QP}.
\end{remark}


\section{EM scheme with exponential family distributions}\label{sec:expEM}
In this section, we show that the EM scheme results in a natural gradient descent algorithm with a fixed step size if we employ an exponential family distribution that satisfies certain conditions. 

\subsection{Basic properties of 
EM scheme with exponential family distributions}
\label{sec:properties.expEM}
A probability distribution of $\bm{X}$ on $\Xi$ 
is said to be (a member of) the exponential family if its probability 
density/mass function, or their hybrid, $p_{\bm{x}}(\bm{\theta})$, of $\bm{x}\in\Xi$ can be rewritten as
\begin{align}
p_{\bm{x}}(\bm{\theta}) &=\beta(\bm{x})\exp\big( {\bm{\gamma}(\bm{x})^\top\bm{\eta}(\bm{\theta})} -A(\bm{\theta})\big)
\label{eq:expdist}
\end{align}
for some functions $A:\mathbb{R}^{p}\to\mathbb{R}$, $\beta:\Xi\to[0,\infty)$ and vectors of functions $\bm{\eta}:\mathbb{R}^{p}\to\mathbb{R}^{q}$, $\bm{\gamma}:\Xi\to\mathbb{R}^{q}$. 
The exponential family encompasses a wide array of representative distributions (see, for instance, \cite{bedbur2021multivariate}). Subsequent sections will demonstrate the significance of certain exponential family distributions, including the normal distribution, binomial distribution, multinomial distribution, and Poisson distribution, in the context of the relationship between EM schemes and polynomial objective optimization problems. First, the basic properties of the exponential family when applied to Algorithm \ref{alg:prototypeEM} are summarized in the following proposition.
\begin{proposition}
\label{prop:char.exp.family}
Suppose that any probability distribution $p_{\bm{x}}(\bm{\theta})$ of the exponential family \eqref{eq:expdist} is used for Algorithm  \ref{alg:prototypeEM} and let $\{\bm{\theta}^{(t)}\}_{t\geq 0}$ denote the sequence of solutions generated by the algorithm. Then we have the following statements.\\  
(a) The probability distribution \eqref{eq:q_x} is represented as
\begin{align}
q_{\bm{x}}(\bm{\theta}^{(t-1)})
 & = \tilde{\beta}(\bm{x}) \exp\left( \bm{\gamma}(\bm{x})^\top\bm{\eta}(\bm{\theta}^{(t-1)}) -\tilde{A}(\bm{\theta}^{(t-1)})\right)
, 
\label{eq:q_exp.family}
\end{align}
where
\begin{align}
 \tilde{A}(\bm{\theta}^{(t-1)}) & :=\ln\mathbb{E}_{p_{\bm{x}}(\bm{\theta}^{(t-1)})}[G(\bm{X})] + A(\bm{\theta}^{(t-1)}), 
 \label{q_exponential_family}
 \\
 \tilde{\beta}(\bm{x}) & := G(\bm{x})\beta(\bm{x}).\nonumber
\end{align}
(b) 
The E-step and M-step (lines 3 to 4) of Algorithm \ref{alg:prototypeEM} can be merged into the minimization of 
the function $u(\bm{\theta} | \bm{\theta}^{(t-1)})$ represented as
\begin{align}
u(\bm{\theta}\vert\bm{\theta}^{(t-1)})
& =  - \mathbb{E}_{q_{\bm{x}}(\bm{\theta}^{(t-1)})}[\bm{\gamma}(\bm{X})]^\top\big(\bm{\eta}(\bm{\theta})-\bm{\eta}(\bm{\theta}^{(t-1)})\big) + A(\bm{\theta}) - A(\bm{\theta}^{(t-1)}) + f(\bm{\theta}^{(t-1)}).
\label{eq:u_component}
\end{align}
(c) Suppose that $A(\bm{\theta})$ and $\bm{\eta}(\bm{\theta})$ are differentiable. If it holds that $\mathbb{E}_{p_{\bm{x}}(\bm{\theta})}[\bm{\gamma}(\bm{X})] =  \mathbb{E}_{q_{\bm{x}}(\bm{\theta}^{(t-1)})}[\bm{\gamma}(\bm{X})]$, then we have 
$\nabla_{\bm{\theta}}u(\bm{\theta}|\bm{\theta}^{(t-1)})=\bm{0}$ where $\bm{0}:=(0,...,0)^\top$.
\\
(d) Suppose that $\bm{\eta}(\bm{\theta})$ is square (i.e., $q=p$) and invertible, and each component $\eta_i(\bm{\theta})$ is differentiable, $i\in[p](=[q])$. 
Then, we have
\begin{align}
\mathbb{E}_{q_{\bm{x}}(\bm{\theta})}\big[\bm{\gamma}(\bm{X}) \big]
& = \mathbb{E}_{p_{\bm{x}}(\bm{\theta})}[\bm{\gamma}(\bm{X})] - \big(\nabla\bm{\eta}(\bm{\theta})^\top \big)^{-1} \nabla f(\bm{\theta}),
\label{eq:Expectation_q}
\end{align}
where 
$\nabla\bm{\eta}(\bm{\theta})$ denotes the Jacobian of $\bm{\eta}(\bm{\theta})$, i.e., $\nabla\bm{\eta}(\bm{\theta}):=(\nabla\eta_1(\bm{\theta}),...,\nabla\eta_p(\bm{\theta}))^\top$.  
\end{proposition}

Statement (a) shows that when the probability distribution $\big(p_{\bm{x}}(\bm{\theta})\big)_{\bm{x}}$ is a member of the exponential family, so is the distribution $\big(q_{\bm{x}}(\bm{\theta})\big)_{\bm{x}}$. Statement (b) states that 
lines 3 and 4 of Algorithm \ref{alg:prototypeEM} can be replaced with 
\begin{align}
\bm{\theta}^{(t)}\gets\arg\min_{\bm{\theta}\in\Theta}\Big\{
A(\bm{\theta})-\mathbb{E}_{q_{\bm{x}}(\bm{\theta}^{(t-1)})}[\bm{\gamma}(\bm{X})]^\top\bm{\eta}(\bm{\theta})
\Big\}.
\label{eq:EM-step}
\end{align}
It also implies that $u(\bm{\theta}\vert\bm{\theta}^{(t-1)})$ is convex in $\bm{\theta}$ and, thus, \eqref{eq:EM-step} is a convex optimization problem if $A(\bm{\theta})-\mathbb{E}_{q_{\bm{x}}(\bm{\theta}^{(t-1)})}[\bm{\gamma}(\bm{X})]^\top\bm{\eta}(\bm{\theta})$ is convex in $\bm{\theta}$. Statement (c) yields a sufficient condition for the stationarity of the M-step \eqref{eq:EM-step}. Statement (d) will be used to derive several formulae, especially, the natural gradient update formula in Theorem \ref{prop:char.exp.family_natgrad}. 

The proof of Proposition \ref{prop:char.exp.family} is given in Appendix \ref{sec:proof.prop:char.exp.family}.  
\begin{corollary}
\label{cor:FOS}
Suppose the same assumption as in Proposition \ref{prop:char.exp.family} and that $f$ and $u$ are differentiable. 
Then, we have 
$\nabla u(\bm{\theta}\vert \bar{\bm{\theta}}) \big|_{\bm{\theta}=\bar{\bm{\theta}}} = \nabla f(\bar{\bm{\theta}})$.
\end{corollary}
A surrogate function $u$ satisfying this property is called first-order surrogate for $f$ (\cite{mairal2013optimization}). 
If $u$ is a first-order surrogate for $f$, it holds true that $u'_{\bm{d}}(\bm{\theta}\vert\bar{\bm{\theta}})\big|_{\bm{\theta}=\bar{\bm{\theta}}}=f'_{\bm{d}}(\bar{\bm{\theta}})$ for all $\bar{\bm{\theta}}\in\Theta$ and $\bm{d}\in\mathbb{R}^p$ such that $\bar{\bm{\theta}}+\bm{d}\in\Theta$, which is a restrictive condition often assumed in the context of statistical EM algorithms, where $f'_{\bm{d}}(\bar{\bm{\theta}})$ denotes the one-sided directional derivative of $f$ at $\bar{\bm{\theta}}$ with respect to the direction $\bm{d}$. 
In addition, this property enables us to apply the Nesterov's acceleration to Algorithm \ref{alg:prototypeEM}. See Appendix \ref{sec:proof:cor:FOS} for the proof.

\subsection{EM scheme with simple exponential family distribution}
Building upon the basic properties presented in the previous subsection, we will show that the update formula 
\eqref{eq:EM-step} for a class of exponential family distributions results in natural gradient descent without line search. 

Let $\bm{I}_{p_{\bm{x}}}(\bm{\theta})$ denote the Fisher information matrix of a probability distribution $p_{\bm{x}}(\bm{\theta})$, i.e.,
\begin{align*}
\bm{I}_{p_{\bm{x}}}(\bm{\theta}):=\mathbb{E}_{p_{\bm{x}}(\bm{\theta})}\Big[ \nabla \ln p_{\bm{X}}(\bm{\theta}) \nabla \ln p_{\bm{X}}(\bm{\theta})^\top \Big].
\end{align*}

\begin{theorem}
\label{prop:char.exp.family_natgrad}
Suppose that 
$p_{\bm{x}}(\bm{\theta})$ is differentiable and 
a member of the exponential family \eqref{eq:expdist} with $\nabla\bm{\eta}(\bm{\theta})$ being a $p\times p$ invertible  matrix, and that there exists an invertible matrix $\bm{\Omega}\in\mathbb{R}^{p\times p}$ such that $\bm{\Omega}\bm{\gamma}(\bm{X})$ is an unbiased estimator of $\bm{\theta}$, 
i.e., 
\begin{align}
 \mathbb{E}_{p_{\bm{x}}(\bm{\theta})}\big[\bm{\Omega}\,\bm{\gamma}(\bm{X})\big]=\bm{\theta}.
 \label{eq:unbiased_estimator}
\end{align}
Then we have 
\begin{align}
\bm{I}_{p_{\bm{x}}}(\bm{\theta}) 
= \big(\mathbb{V}_{p_{\bm{x}}(\bm{\theta})}[\bm{\Omega}\bm{\gamma}(\bm{X})]\big)^{-1}
= \nabla\bm{\eta}(\bm{\theta})^\top \bm{\Omega}^{-1},
\label{fisher_inv}
\end{align}
where $\mathbb{V}_{p_{\bm{x}}(\bm{\theta})}[\cdot]$ denotes the (co)variance operator under the probability distribution $p_{\bm{x}}(\bm{\theta})$ of $\bm{X}$. 
In particular, if $f$ is differentiable, the update formula providing the stationarity point of \eqref{eq:EM-step} 
is then reduced to a natural gradient descent of the form:
\begin{align}
 \bm{\theta}^{(t)}  \gets \bm{\theta}^{(t-1)} - \bm{I}_{p_{\bm{x}}}(\bm{\theta}^{(t-1)})^{-1} \nabla f(\bm{\theta}^{(t-1)}).\label{eq:natgrad}
\end{align}
\end{theorem}
See Appendix \ref{sec:proof:prop:char.exp.family_natgrad} for the proof.  

Intuitively, the natural gradient \eqref{eq:natgrad} modifies the gradient, reflecting the metric of a manifold of the exponential family distribution $p_{\bm{x}}(\bm{\theta})$ via its Fisher information matrix or the covariance matrix of $\bm{\Omega}\bm{\gamma}(\bm{X})$. 
This comes from the fact that we employ $p_{\bm{x}}(\bm{\theta})$ so as to make the general optimization problem \eqref{eq:P} equivalent to minimizing \eqref{eq:G-condition}. 
Next section will present three examples to show how an optimization problem can be connected to a natural gradient descent.

\begin{remark}
Theorem \ref{prop:char.exp.family_natgrad} assumes that $\bm{\Omega}$ does not depend on $\bm{\theta}$. If, on the other hand, $\bm{\Omega}$ depends on $\bm{\theta}$, i.e., $\bm{\Omega}=\bm{\Omega}(\bm{\theta})$, we can see that \eqref{eq:natgrad} becomes somewhat more complicated, such as 
\begin{align*}
 \bm{\theta}^{(t)}  \gets \bm{\Omega}(\bm{\theta}^{(t)})\bm{\Omega}(\bm{\theta}^{(t-1)})^{-1}\big(\bm{\theta}^{(t-1)} - \bm{I}_{p_{\bm{x}}}(\bm{\theta}^{(t-1)})^{-1} \nabla f(\bm{\theta}^{(t-1)})\big).
\end{align*}
All the examples demonstrated below are constructed based on the constant matrix $\bm{\Omega}$, and as a result we can enjoy the simple formula \eqref{eq:natgrad}.
\end{remark}


\section{Three simple examples}\label{sec:examples}
This section presents three applications of the EM scheme under exponential family distributions.


\subsection{Unconstrained convex quadratic minimization}
\label{subsec:unconst.cvx.quadratic}
We start with the unconstrained minimization of a convex quadratic function:
\begin{align}
\begin{array}{ll}
\underset{\bm{\theta}\in\mathbb{R}^{p}}{\mbox{minimize}} \quad & f(\bm{\theta}) = \displaystyle \frac{1}{2}\bm{\theta}^\top\bm{Q}\bm{\theta}+\bm{b}^\top\bm{\theta},
\end{array}
\label{QP}
\end{align}
where $\bm{Q}\in\mathbb{S}^{p}_+$ and $\bm{b}\in\mathbb{R}^p$.  To apply the EM scheme (Algorithm \ref{alg:prototypeEM}) to solve \eqref{QP}, we employ the $p$-dimensional normal distribution $\bm{{\mathrm{N}}}(\bm{\theta},\bm{\Sigma})$ for random variables $\bm{X}=(X_1,...,X_p)$, where $\bm{\Sigma}\in\mathbb{S}^p_{++}$. Here, we suppose that the covariance $\bm{\Sigma}$ is known and the mean $\bm{\theta}$ is the only unknown parameter, which corresponds to the decision variables of the problem \eqref{QP}. The normal distribution is a member of the exponential family \eqref{eq:expdist} since its density can be decomposed as
\begin{align}
p_{\bm{x}}(\bm{\theta})
& = \frac{1}{(2\pi)^{\frac{p}{2}} \sqrt{|\bm{\Sigma}|}} \exp\Big(- \frac{1}{2}(\bm{x}-\bm{\theta})^\top\bm{\Sigma}^{-1} (\bm{x}-\bm{\theta}) \Big)
\label{eq:normal.prob}
\\
& = \underbrace{\exp\big(-\frac{1}{2}\bm{x}^\top\bm{\Sigma}^{-1}\bm{x}\big)}_{=\beta(\bm{x})}
\exp\big(
\underbrace{\bm{x}^\top\bm{\Sigma}^{-1}\bm{\theta}}_{=\bm{\gamma}(\bm{x})^\top\bm{\eta}(\bm{\theta})}
\underbrace{-\ln (2\pi)^{\frac{p}{2}}\sqrt{|\bm{\Sigma}|}-\frac{1}{2}\bm{\theta}^\top\bm{\Sigma}^{-1}\bm{\theta}}_{=-A(\bm{\theta})}
\big).
\label{eq:normal.as.exponential}
\end{align}
\begin{proposition}
\label{prop:G.for.cvx.q.min}
Suppose that $\bm{\Sigma}^{-1}\succ\bm{Q}$, where $\bm{\Sigma}\in\mathbb{S}^p_{++}$ is the covariance matrix of normal distribution \eqref{eq:normal.prob}, and $\bm{Q}\in\mathbb{S}^p_{+}$ is the matrix defining the quadratic objective of \eqref{QP}. 
If $\bm{X}$ follows the normal distribution \eqref{eq:normal.prob} and the functions $G$ is defined by 
\begin{align}
G(\bm{x}) = \sqrt{\frac{|\bm{\Sigma}|}{|\tilde{\bm{Q}} + \bm{\Sigma}^{-1}|}} \exp \left( - \frac{1}{2}\bm{x}^\top \tilde{\bm{Q}}\bm{x} - \tilde{\bm{b}}^\top \bm{x} \right), 
\label{eq:G.for.cvxQP}
\end{align}
where
$
 \tilde{\bm{Q}} :=  \bm{\Sigma}^{-1}(\bm{\Sigma}^{-1}-\bm{Q})^{-1} \bm{\Sigma}^{-1} - \bm{\Sigma}^{-1}
 \quad\mbox{and}\quad
 \tilde{\bm{b}} :=  \bm{\Sigma}^{-1}(\bm{\Sigma}^{-1}-\bm{Q})^{-1} \bm{b},
$ 
the condition \eqref{eq:G-condition} for \eqref{QP} is fulfilled with $p_{\bm{x}}(\bm{\theta})$ defined by \eqref{eq:normal.prob} and $\Theta=\mathbb{R}^p$. 
\end{proposition}
See Appendix \ref{sec:proof:prop:G.for.cvx.q.min} for the proof.  

From the correspondence indicated in \eqref{eq:normal.as.exponential}, $\bm{\gamma}(\bm{x})=\bm{x}$ holds for the normal distribution. 
If we set $\bm{\Omega}=\bm{I}$ for \eqref{eq:unbiased_estimator}, we see $\mathbb{E}_{p_{\bm{x}}(\bm{\theta})}\big[\bm{\Omega} \bm{\gamma}(\bm{X})\big]=\mathbb{E}_{p_{\bm{x}}(\bm{\theta})}[\bm{X}]=\bm{\theta}$, so $\bm{\Omega}\bm{\gamma}(\bm{X})=\bm{\gamma}(\bm{X})$ turns out to be an unbiased estimator of $\bm{X}$. Therefore, the update formula \eqref{eq:natgrad} for \eqref{eq:EM-step} is simplified as follows.
\begin{corollary} \label{cor:updateQP}
If the normal distribution \eqref{eq:normal.prob} is employed for Algorithm \ref{alg:prototypeEM} to solve \eqref{QP}, the update formula \eqref{eq:natgrad} results in
\begin{align} 
\bm{\theta}^{(t)} \gets \bm{\theta}^{(t-1)} - \bm{\Sigma}
 (\bm{Q}\bm{\theta}^{(t-1)}+\bm{b}).
\label{eq:updateUnconstQP}
\end{align}
\end{corollary}
Algorithm \ref{alg:EMforUnconstQP} shows the full description of the procedure. 
\begin{algorithm}[tb]
\caption{EM algorithm for unconstrained convex quadratic minimization \eqref{QP}}
\label{alg:EMforUnconstQP}
\begin{algorithmic}[1]
\State 
Let $\bm{\Sigma}\in\mathbb{S}^p_{++}$ satisfy $\bm{\Sigma}^{-1}\succ\bm{Q}$, and $\bm{\theta}^{(0)}\in\mathbb{R}^n$, and set $t\gets 1$.
\Repeat
\State
Update the solution by \eqref{eq:updateUnconstQP} and $t\gets t+1$. 
\Until{a termination condition is fulfilled}
\end{algorithmic}
\end{algorithm}
From the MM optimization perspective shown in Section \ref{sec:framework}
, the natural gradient formula \eqref{eq:updateUnconstQP} comes from the minimization of the convex quadratic surrogate function 
\begin{align}
u(\bm{\theta}\vert\bm{\theta}^{(t-1)})
 &= \frac{1}{2}\bm{\theta}^\top\bm{\Sigma}^{-1}\bm{\theta} 
 - \big\{ (\bm{\Sigma}^{-1}-\bm{Q})\bm{\theta}^{(t-1)}-\bm{b}
 \big\}^\top \big(\bm{\theta}-\bm{\theta}^{(t-1)}\big)\nonumber\\
 &\quad\qquad - \frac{1}{2}(\bm{\theta}^{(t-1)})^\top\bm{\Sigma}^{-1}\bm{\theta}^{(t-1)} + f(\bm{\theta}^{(t-1)}),
\label{eq:quadratic_surrogate}
\end{align}

which can be derived from Proposition \ref{prop:char.exp.family} (b) and whose convexity is stronger than $f(\bm{\theta})$ because of the relation of the Hessian matrices, $\nabla^2 u(\bm{\theta}|\bar{\bm{\theta}})=\bm{\Sigma}^{-1}\succ\bm{Q}=\nabla^2f(\bm{\theta})$. 
The quadratic surrogate function \eqref{eq:quadratic_surrogate}, although it was not originally intended to be applied to quadratic optimization in that context, can be considered as a special case of the quadratic lower bound principle \cite{bohning1988monotonicity}. 
In addition, observing that \eqref{eq:quadratic_surrogate} is decomposed as
\begin{align*}
u(\bm{\theta}\vert\bm{\theta}^{(t-1)})
 = f(\bm{\theta})+
 \frac{1}{2}\|\bm{\theta}-\bm{\theta}^{(t-1)}\|^2_{\bm{\Sigma}^{-1}-\bm{Q}},
\end{align*}
the MM algorithm can be viewed as a proximal point algorithm with an (elliptical) quadratic proximal term. 
Recalling that the surrogate functions in our EM scheme are all derived from concavity of the logarithm as in \eqref{eq:how_to_majorize}, it may be interesting to note that different majorization principles lead to the same surrogate function. 

To see the role of the covariance matrix $\bm{\Sigma}$ in the formula \eqref{eq:updateUnconstQP}, assume that  
$\bm{Q}\in\mathbb{S}^p_{++}$ 
and let $\bm{\Sigma}=\bm{Q}^{-1}$, hypothetically. This case does not satisfy the assumption of Proposition \ref{prop:G.for.cvx.q.min}, but, for such $\bm{\Sigma}$, we could get 
the update formula 
as $\bm{\theta}^{(t)}\gets \bm{\theta}^{(t-1)} - \bm{Q}^{-1}(\bm{Q}\bm{\theta}^{(t-1)}+\bm{b})= - \bm{Q}^{-1}\bm{b}$, which is equivalent to the solution of the problem \eqref{QP}. 
For a general $\bm{\Sigma}\in\mathbb{S}^p_{++}$, the formula \eqref{eq:updateUnconstQP} can be viewed as a generalization of the steepest descent, i.e., $\bm{\Sigma}=\sigma\bm{I}$ for a constant $\sigma>0$. 
As a bit more practical example is a diagonal matrix $\bm{\Sigma}=\mathrm{diag}(\sigma_1,...,\sigma_p)$ for $\sigma_j>0,j\in[p]$. Note that, for example, setting $\sigma_{j}\in(0,1/(q_{jj}+\sum_{h\neq j}
|q_{jh}|)),j\in[p],$ makes $\bm{\Sigma}^{-1}-\bm{Q}$ diagonally dominant, thus ensuring positive definite (e.g., \cite{horn2012matrix}). Despite that it permits the step size to vary for each coordinate, the practical advantage of this modification over the steepest descent may not be as significant as one might expect. However, the relationship between the convex quadratic objective and the normal distribution described above can be used for the example of Section \ref{sec:l1QP} below. 
Note also that while we here assume $\bm{Q}\in\mathbb{S}^p_+$ so that \eqref{QP} will be well-defined, the convexity of the quadratic objective is not crucial for the MM scheme. In fact, in Section \ref{subsec:QP} we will consider a possibly nonconvex quadratic objective in developing the EM algorithm for constrained problems with the help of the normal distribution.

Convergence rate of Algorithm \ref{alg:EMforUnconstQP} can be shown by applying existing results on EM algorithms \cite{kumar2017convergence,karimi2016linear}.
\begin{corollary}
Let $\overline{\lambda}_{\bm{S}},\underline{\lambda}_{\bm{S}}$ denote the largest and smallest eigenvalues of a real symmetric matrix $\bm{S}$, respectively. 
Suppose the same assumptions as in Corollary \ref{cor:updateQP}. 
(i) We have 
$
 f(\bm{\theta}^{(t-1)}) - f(\bm{\theta}^{(t)}) \geq 
 \| \bm{\theta}^{(t)} - \bm{\theta}^{(t-1)} \|^2
 /(2\overline{\lambda}_{\bm{\Sigma}}),
$
which implies 
$
 \min_{t \in [T]
  } \| \bm{\theta}^{(t)} - \bm{\theta}^{(t-1)} \|^2 \leq
2\overline{\lambda}_{\bm{\Sigma}}\big(f(\bm{\theta}^{(0)}) - f(\bm{\theta}^*)\big)/T.
$ 
(ii) We have 
$ 
f(\bm{\theta}^{(t-1)}) - f(\bm{\theta}^{(t)}) \geq 
(\underline{\lambda}_{\bm{\Sigma}}/2)
\| \nabla f(\bm{\theta}^{(t-1)}) \|^2,
$
which implies 
$
 \min_{t \in [T]} \| \nabla f(\bm{\theta}^{(t-1)}) \|^2 \leq 
(2/\underline{\lambda}_{\bm{\Sigma}})\big(f(\bm{\theta}^{(0)}) - f(\bm{\theta}^*)\big)/T.
$
(iii) In addition, if $\underline{\lambda}_{\bm{Q}}>0$, we have 
$
f(\bm{\theta}^{(T)}) - f(\bm{\theta}^*) \leq \left( 1 - \underline{\lambda}_{\bm{Q}}\underline{\lambda}_{\bm{\Sigma}} \right)^T \big( f(\bm{\theta}^{(0)}) - f(\bm{\theta}^*) \big).
$
\end{corollary}
The statement (i) follows from Theorem 1 of \cite{kumar2017convergence}.  
The proof of (ii) goes similar to that of Theorem 1 of \cite{karimi2016linear} with additional components such as the function $G$. See Theorem 1 of \cite{karimi2016linear} as well for the statement (iii). In addition, we can adopt the Nesterov's acceleration technique to Algorithm \ref{alg:prototypeEM} because this case satisfies the first-order surrogate condition 
(see Corollary \ref{cor:FOS} in Section \ref{sec:properties.expEM}). 

\subsubsection{Extension: $\ell_1$-regularized convex quadratic minimization}
\label{sec:l1QP}
With $f(\bm{\theta})=\frac{1}{2}\bm{\theta}^\top\bm{Q}\bm{\theta}+\bm{b}^\top\bm{\theta}$ and $\rho(\bm{\theta})=\|\bm{\theta}\|_1:=\sum_{j\in[p]}|\theta_j|$, the composite optimization \eqref{eq:F} becomes the $\ell_1$-regularized convex quadratic minimization:
\begin{align*}
\underset{\bm{\theta}\in\mathbb{R}^{p}}{\mbox{minimize}}~
\Phi(\bm{\theta}) = \frac{1}{2}\bm{\theta}^\top\bm{Q}\bm{\theta}+\bm{b}^\top\bm{\theta}
+ \|\bm{\theta}\|_1.
\end{align*}
With the diagonal matrix $\bm{\Sigma}=\mathrm{diag}(\sigma_1,...,\sigma_p)$, the surrogate function becomes separable as 
$U(\bm{\theta}\vert\bar{\bm{\theta}}) =  u(\bm{\theta}\vert\bar{\bm{\theta}}) + \|\bm{\theta}\|_1 
=  \sum_{j\in[p]} \big( 
(\theta_j - \mathbb{E}_{q_{\bm{x}}(\bar{\bm{\theta}})}[X_j])^2
/(2\sigma_j)
+ |\theta_j| \big) + \mbox{``constant,''} $ and it can be minimized via the so-called soft-thresholding operation, i.e., for 
$
\mathbb{E}_{q_{\bm{x}}(\bar{\bm{\theta}})}[X_j] = \bar{\theta}_j - \sigma_j \sum_{h\in[p]}(q_{jh}\bar{\theta}_h+b_j),~j\in [p],
$ 
the coordinate-wise update formula is derived as $\theta_j^{(t)} \gets
{\rm soft}_{\sigma_j}(\mathbb{E}_{q_{\bm{x}}(\bm{\theta}^{(t-1)})}[X_j])$ where ${\rm soft}_{a}(z) := \arg\min_x \left\{ a\vert x \vert + \frac{1}{2}(x-z)^2 \right\} =
z + a$ if $z \leq - a$;~$0$ if $-a \leq z \leq a$;~$z - a$ if $a \leq z$, defines the soft-thresholding with parameter $a\in\mathbb{R}$. 
This procedure is analogous to ISTA (Iterative Soft-Thresholding Algorithm) for the $\ell_1$-regularized convex quadratic minimization, yet it differs from the 
vanilla ISTA in that it permits the step size to vary for each coordinate. Also, while ISTA has famous applications in image restoration and reconstruction, and EM algorithms associated with normal distribution exist (e.g., \cite{figueiredo2003algorithm}), the role of normal distribution in the above treatment has little to do with that in such statistical contexts. 

\subsection{Minimization of a polynomial over a rectangle}
\label{sec:poly.on.rect}
In this subsection, we consider the minimization of a polynomial function 
over 
a hyper-rectangle in $\mathbb{R}^{p}$: 
\begin{align}
\begin{array}{|ll}
\underset{\bm{\lambda}\in\mathbb{R}^{p}}{\mbox{minimize}} \quad & F(\bm{\lambda}):=\sum\limits_{i\in[I]}
a_i \prod\limits_{j\in[p]} \lambda_{j}^{n_{i,j}} \\
\mbox{subject to} & \bm{l}\leq\bm{\lambda}\leq\bm{u},
\end{array}
\label{box_multi_problem}
\end{align}
where $a_i\in\mathbb{R}$ for $i\in [I]$, $n_{i,j}\in\{0\}\cup
\mathbb{N}$ for $i\in [I],j\in[p]$, and $\bm{l},\bm{u}\in\mathbb{R}^p$ satisfy $l_j<u_j$ for each $j$. 
For example, the dual form of Support Vector Classification for a bias-free decision function can be cast in this form with a convex quadratic objective. 
To establish a binomial distribution-based EM algorithm for this optimization problem, we first note that by the change of variables $\theta_j:=(\lambda_j-l_j)/(u_j-l_j)$, the problem \eqref{box_multi_problem} can be rewritten as
\begin{align}
\begin{array}{|ll}
\underset{\bm{\theta}\in\mathbb{R}^{p}}{\mbox{minimize}} \quad & \sum\limits_{i\in[I]} 
 \tilde{a}_i \prod\limits_{j\in[p]} \theta_j^{n_{i,j}} + \tilde{b}\quad(=F\big(\mathrm{diag}(\bm{u}-\bm{l})\bm{\theta}+\bm{l}\big)) \\
\mbox{subject to} & \bm{0}\leq\bm{\theta}\leq\bm{1},
\end{array}
\label{box_multi_problem_theta}
\end{align}
for some $\tilde{a}_i\in\mathbb{R},i\in[I]$, and $\tilde{b}\in\mathbb{R}$ such that
\begin{align*}
 F(\mathrm{diag}(\bm{u}-\bm{l})\bm{\theta}+\bm{l}) & = \sum_{i\in[I]}  a_i \prod_{j\in[p]} \big((u_j-l_j)\theta_j+l_j\big)^{n_{i,j}} = \sum_{i\in[I]} \tilde{a}_i \prod_{j\in[p]} \theta_j^{n_{i,j}} + \tilde{b}.
\end{align*}
Note also that there exists an upper bound $K$ of the objective value of \eqref{box_multi_problem} and \eqref{box_multi_problem_theta} because the objective function is continuous and the feasible region is compact. 
With such an upper bound $K$, i.e., $K>\max\{F(\bm{\lambda})|\bm{\lambda}\in[\bm{l},\bm{u}]\}(=\max\{F(\mathrm{diag}(\bm{u}-\bm{l})\bm{\theta}+\bm{l})|\bm{\theta}\in[\bm{0},\bm{1}]\})$, 
\eqref{box_multi_problem_theta} 
can be equivalently rewritten as 
\begin{align}
\begin{array}{|ll}
\underset{\bm{\theta}\in\mathbb{R}^p}{\mbox{minimize}}\quad
 & f(\bm{\theta}) = -\ln \Big(K-F\big(\mathrm{diag}(\bm{u}-\bm{l})\bm{\theta}+\bm{l}\big)\Big) \\
\mbox{subject to}
 & \bm{0}\leq\bm{\theta}\leq\bm{1}.
\end{array}
\label{box_multi_problem2}
\end{align}
See Appendix \ref{sec:Setting.of.K.for.rectangle} for an example of how to obtain such $K$. 

To apply the EM scheme to \eqref{box_multi_problem2}, 
we consider a vector $\bm{X}=(X_1,...,X_p)$ 
of $p$ independent binomial random variables, $X_j\sim\mathrm{Bin}(m_j,\theta_j)$, with $m_j\geq\max\{n_{i,j}\,|\,i\in
[I]\},~j\in [p]$, i.e., $m_j$ must be no less than the highest degree of variable $\theta_j$. $\bm{X}$ is then of the exponential family \eqref{eq:expdist} since its probability mass on $\bm{x}\in\Xi=\prod_{j\in[p]}(\{0\}\cup[m_j])$ is given by 
\begin{align}
p_{\bm{x}}(\bm{\theta})
 &=\prod_{j\in[p]}
 \frac{m_j!}{x_j! (m_j-x_j)!}\theta_j^{x_j}(1-\theta_j)^{m_j-x_j} \label{eq:binomial.prob}
\\
 & =\underbrace{\prod_{j\in[p]}\frac{m_j!}{x_j! (m_j-x_j)!}}_{=\beta(\bm{x})}
   \exp\Big(\underbrace{\sum_{j\in[p]} x_j\ln\frac{\theta_j}{1-\theta_j}}_{=\bm{\gamma}(\bm{x})^\top\bm{\eta}(\bm{\theta})}+\underbrace{\sum_{j\in[p]}m_j\ln(1-\theta_j)}_{=-A(\bm{\theta})}\Big).\nonumber
\end{align}
The next proposition shows a function $G(\cdot)$ that satisfies \eqref{eq:G-condition} for \eqref{box_multi_problem2} under \eqref{eq:binomial.prob}. 
\begin{proposition}
\label{prop:expectation_condition_for_simple_binom}
Let $m_j\geq\max\{n_{i,j}\,|\,i\in
[I]\}$ for $j\in [p]$. 
If $G$ is defined by
\begin{align}
G(\bm{x}) & = K - \Big( \sum_{i\in[I]} \tilde{a}_i \prod_{j\in[p]}
g_{i,j}(x_j) + \tilde{b} \Big),
\label{def:G}
\end{align}
where, corresponding to $n_{i,j}$, $g_{i,j}(\cdot)$ is defined with $m_j$ 
as
\begin{align*}
 g_{i,j}(x_j) := &
\begin{cases}
\displaystyle \frac{x_j(x_j-1) \cdots (x_j-n_{i,j}+1)}{m_j(m_j-1) \cdots (m_j-n_{i,j}+1)} , & \mbox{if }n_{i,j}>0, \\
 1, & \mbox{if }n_{i,j} =0,
\end{cases}
\end{align*}
the condition \eqref{eq:G-condition} for \eqref{box_multi_problem2} is fulfilled with $p_{\bm{x}}(\bm{\theta})$ defined by \eqref{eq:binomial.prob} and $\Theta=[\bm{0},\bm{1}]$.
\end{proposition}
See Appendix \ref{sec:proof:prop:expectation_condition_for_simple_binom} for the proof. 

Applying Proposition \ref{prop:char.exp.family} (b), we see that 
the EM algorithm iteratively minimizes the surrogate function 
\begin{align*}
\lefteqn{u(\bm{\theta}\vert\bm{\theta}^{(t-1)})}\\
&= 
 - \sum\limits_{j\in[p]}
\Big\{
\mathbb{E}_{q_{\bm{x}}(\theta^{(t-1)}_j)}[X_j] \ln\frac{\theta_j}{\theta^{(t-1)}_j}
+
(m_j-\mathbb{E}_{q_{\bm{x}}(\theta^{(t-1)}_j)}[X_j]) \ln\frac{1-\theta_j}{1-\theta^{(t-1)}_j}
\Big\}
+ f(\bm{\theta}^{(t-1)}),
\end{align*}
which can be reduced to a natural gradient descent as follows.
\begin{corollary}
\label{cor:update_formula_binomial}
If the binomial distribution \eqref{eq:binomial.prob}  
is employed for Algorithm \ref{alg:prototypeEM} 
to solve \eqref{box_multi_problem2}, the update formula \eqref{eq:natgrad} 
results in
\begin{align}
\theta^{(t)}_j\gets\theta^{(t-1)}_j-\frac{\theta^{(t-1)}_{j}(1-\theta^{(t-1)}_{j})}{m_j}\frac{\partial }{\partial\theta_j}f(\bm{\theta}^{(t-1)}),\quad j\in[p].
\label{eq:update.formula.poly.on.rect}
\end{align}
\end{corollary}
See Appendix \ref{sec:proof_update_formula_binomial} for the proof.

Corollary \ref{cor:update_formula_binomial} reveals that minimizing the surrogate function (or \eqref{eq:EM-step}) for \eqref{box_multi_problem2} results in a component-wise gradient descent \eqref{eq:update.formula.poly.on.rect} whose step size is endogenously determined at each iteration as $\theta_{j}^{(t-1)}\big(1-\theta_{j}^{(t-1)}\big)/m_{j},~j \in [p]$. 
Given that $m_j$ may take any number satisfying $m_j\geq\max\{n_{i,j}|i\in[I]\}$, the use of $m_j=\max\{n_{i,j}|i\in[I]\}$ yields the largest step size. 
Note also that if we start the algorithm from an interior point of the rectangle, the generated sequence $\{\bm{\theta}^{(t)}\}_{t\geq 0}$ will remain in the interior, $(0,1)^p$, which implies that the above algorithm is an interior point method. Algorithm \ref{alg:polynom.on.rect} shows the full description of the procedure. 
\begin{algorithm}[tb]
\caption{EM algorithm for minimizing a polynomial objective over a rectangle \eqref{box_multi_problem2}}
\label{alg:polynom.on.rect}
\begin{algorithmic}[1]
\State Let $K >\max_{\bm{\lambda}\in [\bm{l},\bm{u}]}F(\bm{\lambda})$ and $\bm{\theta}^{(0)}\in(0,1)^p$
, and set $t\gets 1$.
\Repeat
\State
Update the solution by \eqref{eq:update.formula.poly.on.rect} and $t\gets t+1$. 
\Until{a termination condition is fulfilled}
\end{algorithmic}
\end{algorithm}


\subsection{Minimization of a polynomial over a unit simplex}
\label{sec:poly.on.simplex}
As the third example, this section considers the minimization of 
the same objective function as 
\eqref{box_multi_problem} 
over a unit simplex. 
\begin{align}
\begin{array}{|ll}
\underset{\bm{\lambda}\in\mathbb{R}^{p}}{\mbox{minimize}} \quad & F(\bm{\lambda}):=\displaystyle \sum\limits_{i\in[I]} a_i \prod\limits_{j\in[p]} \lambda_{j}^{n_{i,j}} \\
\mbox{subject to}
 & \bm{1}^\top\bm{\lambda}=1, ~\bm{\lambda}\geq\bm{0}.\\
\end{array}
\label{min_var_problem}
\end{align}
\eqref{min_var_problem} encompasses, for example, a portfolio selection problem in which the objective function consisting of the weighted sum of higher order moments of asset returns is optimized without a short sale. 
To establish a multinomial distribution-based EM algorithm for \eqref{min_var_problem}, deleting the $p$-th variable $\lambda_p$ by using the relation $\lambda_p=1-\sum_{j\in[p-1]}\lambda_j$ and replacing $\lambda_j$ with $\theta_j$, $j\in[p-1]$, 
we get the equivalent formulation as
\begin{align}
\begin{array}{|ll}
\underset{\theta_1,...,\theta_{p-1}}{\mbox{minimize}} \quad & \displaystyle \sum\limits_{i\in[I]} a_i \prod\limits_{j\in[p-1]} \theta_{j}^{n_{i,j}}\big(1-\sum_{h\in[p-1]}\theta_h\big)^{n_{i,p}} \\
\mbox{subject to}
 & \displaystyle\sum\limits_{j\in[p-1]}\theta_j\leq 1,\quad\theta_j\geq 0,~j\in[p-1].\\
\end{array}
\label{min_var_problem_reduc}
\end{align}
Similarly to \eqref{box_multi_problem2}, 
with an upper bound $K$ of the objective value, the problem \eqref{min_var_problem_reduc} can be equivalently rewritten as
\begin{align}
\begin{array}{|ll}
\underset{\bm{\theta}\in\mathbb{R}^{p-1}}{\mbox{minimize}}\quad
 & f(\bm{\theta}) = -\ln\big(K-F((\bm{I},-\bm{1})^\top\bm{\theta}+\bm{e}_{p})\big) \\
\mbox{subject to}
 & \bm{1}^\top\bm{\theta}\leq 1, ~\bm{\theta}\geq\bm{0},
\end{array}
\label{min_var_problem2}
\end{align}
where the size of $\bm{I}$ and $\bm{1}$ 
is 
$p-1$, i.e., $\bm{I}\in\mathbb{R}^{(p-1)\times(p-1)},\bm{1}\in\mathbb{R}^{p-1}$, and $\bm{e}_{p}$ denotes the unit vector whose $p$-th component is 1, i.e., $\bm{e}_{p}=(0,\ldots,0,1)^\top\in\mathbb{R}^p$. 
See Appendix \ref{sec:Setting.of.K.for.simplex} for an example how to get $K$. 

To apply the EM scheme to \eqref{min_var_problem2}, 
consider the $(p-1)$-dimensional random variables 
$\bm{X} = (X_{1},\ldots,X_{p-1})$ that follows a multinomial distribution, $\bm{X}\sim\mathrm{Multi}(m,\bm{\theta})$, with $m\geq\max\{\sum_{j\in[p]}n_{i,j}\vert i\in[I]\}$. It is of the exponential family since 
its mass on $\bm{x}\in\Xi= \{ \bm{x} \in (\{0\} \cup \mathbb{N})^{p-1} \,|\, \bm{1}^\top\bm{x} \leq m\}$ is given by 
\begin{align}
p_{\bm{x}}(\bm{\theta})
 &=\frac{m!}{
 x_{1}! \cdots x_{p-1}!\big(m-\bm{1}^\top\bm{x}
 \big)!} \prod_{j\in[p-1]}\theta_{j}^{x_{j}} {\big(1- \bm{1}^\top\bm{\theta}
 \big)^{m-\bm{1}^\top\bm{x}
 }
 }
 \label{eq:multinomial.prob} \\
 &=\underbrace{\frac{m!}{x_{1}! \cdots x_{p-1}! \big(m-\bm{1}^\top\bm{x}\big)!}}_{=\beta(\bm{x})} \exp\Big(\underbrace{\sum_{j\in[p-1]}x_j\ln\frac{\theta_j}{1-\bm{1}^\top\bm{\theta}}}_{=\bm{\gamma}(\bm{x})^\top\bm{\eta}(\bm{\theta})}
 \underbrace{+m\ln\big(1-\bm{1}^\top\bm{\theta}\big)}_{=-A(\bm{\theta})}
 \Big).\nonumber
\end{align}
The next proposition shows a function $G(\cdot)$ that satisfies \eqref{eq:G-condition} for \eqref{min_var_problem2} under \eqref{eq:multinomial.prob}.
\begin{proposition}\label{prop:EMcond_multinomial}
Let $m\in\mathbb{N}$ satisfy $m\geq \max\{\sum_{j\in[p]}n_{i,j}|i\in[I]\}$. 
If $G$ is defined by
\begin{align*}
G(\bm{x}) & = K - \sum_{i\in[I]} a_i g_{i}(\bm{x}), 
\end{align*}
where $g_{i}(\cdot)$ is defined with $m$ as 
\begin{align*}
g_{i}(\bm{x}) := &
\left\{\begin{array}{ll}
 \displaystyle \frac{\prod\limits_{j\in[p]}
 x_j(x_j-1) \cdots (x_j-n_{i,j}+1)}{
 m(m-1) \cdots\Big(m-\sum\limits_{j\in[p]} n_{i,j}+1\Big)
 }, & \mbox{if }\sum\limits_{j\in[p]}n_{i,j}> 0,\\
 1, & \mbox{if }\sum\limits_{j\in[p]}n_{i,j} =0,
\end{array}\right.
\end{align*}
where $\bm{x}=(x_1,\ldots,x_{p-1})^\top$ and $x_p = m-\bm{1}^\top\bm{x}$, 
the condition \eqref{eq:G-condition} for \eqref{min_var_problem2} is fulfilled with $ p_{\bm{x}}(\bm{\theta})$ defined by \eqref{eq:multinomial.prob} and $\Theta=\{\bm{\theta}\in\mathbb{R}_+^{p-1}\vert\bm{1}^\top\bm{\theta}\leq 1\}$.
\end{proposition}
See Appendix \ref{sec:proof:prop:expectation_condition_for_simple_multinom} for the proof.

Similarly to the optimization over a rectangle discussed in the previous subsection, the surrogate function is given as
\begin{align*}
\lefteqn{u(\bm{\theta}\vert\bm{\theta}^{(t-1)})= - \sum\limits_{j\in[p-1]} 
\Big\{
\mathbb{E}_{q_{\bm{x}}(\theta^{(t-1)}_j)}[X_j] \ln\frac{\theta_j}{\theta^{(t-1)}_j}\Big\}
}\\
&\qquad\qquad\qquad\quad
-
\Big(m - \sum\limits_{j\in[p-1]}\mathbb{E}_{q_{\bm{x}}(\theta^{(t-1)}_j)}[X_j]\Big) \ln\frac{1-
\bm{1}^\top\bm{\theta}}{1-
\bm{1}^\top\bm{\theta}^{(t-1)}} + f(\bm{\theta}^{(t-1)}),
\end{align*}
and its minimization \eqref{eq:EM-step} can be reduced to a natural gradient descent \eqref{eq:natgrad} as follows.
\begin{corollary}
\label{cor:nat.grad.formula.for.polynom.simplex}
If the multinomial distribution \eqref{eq:multinomial.prob} is employed for Algorithm \ref{alg:prototypeEM} to solve \eqref{min_var_problem2}, the update formula \eqref{eq:natgrad} results in 
\begin{align}
\begin{split}
V^{(t-1)}&\gets\sum_{j\in[p-1]}\theta_j^{(t-1)}\frac{\partial}{\partial\theta_j} f(\bm{\theta}^{(t-1)}),\\
\theta_j^{(t)}&\gets\theta_j^{(t-1)}-\frac{\theta^{(t-1)}_j}{m}\Big(\frac{\partial}{\partial\theta_j} f(\bm{\theta}^{(t-1)})-V^{(t-1)}\Big),\quad j\in[p-1].
\end{split}
\label{eq:update.formula.poly.on.unit}
\end{align}
\end{corollary}
Similarly to \eqref{eq:update.formula.poly.on.rect}, the update formula \eqref{eq:update.formula.poly.on.unit} is component-wise; no line search is involved  
and the use of $m=\max\{\sum_{j\in[p]}n_{i,j}|i\in[I]\}$ yields the largest step size. 
See Appendix \ref{sec:proof:cor:nat.grad.formula.for.polynom.simplex} for the proof.  Algorithm \ref{alg:polynom.on.simplex} shows the full description of the procedure. 
\begin{algorithm}[tb]
\caption{EM algorithm for minimizing the polynomial objective over a simplex \eqref{min_var_problem2}}
\label{alg:polynom.on.simplex}
\begin{algorithmic}[1]
\State Let $K>\max\{F(\bm{\lambda})|\bm{1}^\top\bm{\lambda}=1,\bm{\lambda}\geq\bm{0}\}$ and $\bm{\theta}^{(0)}\in\{\bm{\theta}\in\mathbb{R}_{++}^{p-1}\vert\bm{1}^\top\bm{\theta}< 1\}$, and set $t\gets 1$.
\Repeat
\State
Update the solution by \eqref{eq:update.formula.poly.on.unit} 
and $t\gets t+1$. 
\Until{a termination condition is fulfilled}
\end{algorithmic}
\end{algorithm}


\subsection{Global convergence}\label{sec:convergence}
Global convergence of Algorithms \ref{alg:polynom.on.rect} and \ref{alg:polynom.on.simplex} can be shown in a unified manner. 
Assume that the feasible set of \eqref{eq:P} is represented by a system of inequalities as $\Theta=\{\bm{\theta}\in\mathbb{R}^p\vert g_i(\bm{\theta})\leq 0, i\in[M]\}$ and has an interior point, i.e., $\exists\bm{\theta}$ such that $g_i(\bm{\theta})<0, i\in[M]$:
\begin{align}
\begin{array}{|lll}
\underset{\bm{\theta}\in\mathbb{R}^p}{\mbox{minimize}}\quad& f(\bm{\theta}) \\
\mbox{subject to}~~& g_i(\bm{\theta})\leq0,\quad i\in[M],
\end{array}
\label{eq:shift-form2_theta}
\end{align}
where $g_i:\mathbb{R}^p\to\mathbb{R},i\in[M],$ are continuous. 
\begin{theorem} \label{theorem:global_convergence}
Suppose that $p_{\bm{x}}(\bm{\theta})$ is a member of the exponential family \eqref{eq:expdist} with 
\begin{align}
\bm{\eta}(\bm{\theta}) = \Big(\sum\limits_{i\in[M]} 
 r_{1,i}\ln(-g_i(\bm{\theta})),\ldots,\sum\limits_{i\in[M]} 
 r_{p,i}\ln(-g_i(\bm{\theta})) \Big)^\top \in \mathbb{R}^{p} \label{eq:def_eta_ineq}
\end{align}
for some constants $r_{j,i} \in \mathbb{R},i\in[M];j\in[p]$. 
Suppose that a sequence of solutions, $\{\bm{\theta}^{(t)}\}_{t\geq 0}$, generated by \eqref{eq:EM-step} has the limit $\bm{\theta}^*$, and that for each $i\in [M]$, 
\begin{align}
\lim_{t\to\infty}\frac{\sum_{j\in[p]}r_{j,i}\big(\mathbb{E}_{q_{\bm{x}}(\bm{\theta}^{(t)})}[\bm{\gamma}(\bm{X})]-\mathbb{E}_{q_{\bm{x}}(\bm{\theta}^{(t-1)})}[\bm{\gamma}(\bm{X})]\big)_j}{g_i(\bm{\theta}^{(t)})} \geq 0. 
\label{eq:ineq_limit}
\end{align}
Then $\bm{\theta}^*$ satisfies the Karush--Kuhn--Tucker condition for \eqref{eq:shift-form2_theta}. 
\end{theorem}
See Appendix \ref{sec:proof.Thm:convergence} for the proof. Both Algorithms \ref{alg:polynom.on.rect} and \ref{alg:polynom.on.simplex} satisfy the condition of Theorem \ref{theorem:global_convergence}. 
See Appendix \ref{sec:fulfillment_of_convergence_cond} for the details.  


\section{Extended applications}\label{sec:further_examples}
In this section, we consider extending the EM scheme for a more general optimization formulation of those examined in the previous section. Specifically, for  
a matrix $\bm{H}\in\mathbb{R}^{p\times n}$, a vector $\bm{w}\in\mathbb{R}^p$, and a function $\phi$ on $\mathbb{R}^p$, 
we consider 
\begin{align}
\begin{array}{|lll}
\underset{\bm{\theta}}{\mbox{minimize}}& f(\bm{\theta})=\phi(\bm{H}\bm{\theta}+\bm{w})\\
\mbox{subject to}& \bm{H}\bm{\theta}+\bm{w}\in\Theta,
\end{array}
~\mbox{or
,}~~
\begin{array}{|lll}
\underset{\bm{\theta},\bm{\lambda}}{\mbox{minimize}}& \phi(\bm{\lambda})~\big(=f(\bm{\theta})\big)\\
\mbox{subject to}& \bm{\lambda}=\bm{H}\bm{\theta}+\bm{w},\bm{\lambda}\in\Theta.
\end{array}
\label{eq:shift-form2}
\end{align}
Two formulations in \eqref{eq:shift-form2} are equivalent; for 
$\bm{H}=\bm{I}$ and $\bm{w}=\bm{0}$, they are both reduced to \eqref{eq:P}, and 
each of 
\eqref{QP}, \eqref{box_multi_problem2}, and \eqref{min_var_problem2} can be 
cast as a special case of this extension. 
However, to extend the EM scheme to problems of this general type \eqref{eq:shift-form2}, condition \eqref{eq:G-condition} needs to be modified. 
For \eqref{eq:shift-form2}, consider a distribution $\hat{p}_{\bm{x}}(\bm{\theta}):=p_{\bm{x}}(\bm{H}\bm{\theta}+\bm{w})$ that satisfies $\mathrm{int}\{\bm{\theta}\in\mathbb{R}^n| \bm{H}\bm{\theta}+\bm{w}\in\Theta\}=\{\bm{\theta}\in\mathbb{R}^n\,\vert\,\hat{p}_{\bm{x}}(\bm{\theta})>0\}$. 
 For such $\hat{p}_{\bm{x}}$, assume that there is a positive-valued function $G:\Xi\to(0,\infty)$ satisfying 
\begin{align}
f(\bm{\theta})=-\ln
\mathbb{E}_{\hat{p}_{\bm{x}}(\bm{\theta})}\big[G(\bm{X})
\big]\quad\mbox{for all }\bm{\theta}\in\mathbb{R}^n\mbox{ such that }
\hat{p}_{\bm{x}}(\bm{\theta})
> 0.
\label{eq:G-condition_ext}
\end{align}
If we find such a tuple $(f,\Theta,\hat{p}_{\bm{x}},G)$, the EM scheme can be extended because $\hat{p}_{\bm{x}}(\bm{\theta})=p_{\bm{x}}(\bm{H}\bm{\theta}+\bm{w})=\beta(\bm{x})\exp\big( \bm{\gamma}(\bm{x})^\top\bm{\eta}(\bm{H} \bm{\theta} + \bm{w})-A(\bm{H} \bm{\theta} + \bm{w})\big)$ is a member of the exponential family, and the Fisher information matrix of $\hat{p}_{\bm{x}}(\bm{\theta})$ (or $p_{\bm{x}}(\bm{H} \bm{\theta}+\bm{w})$) is given by 
$
\bm{I}_{\hat{p}_{\bm{x}}}(\bm{\theta})=\bm{H}^\top \big(\mathbb{V}_{p_{\bm{x}}(\bm{H}\bm{\theta}+\bm{w})}[\bm{\Omega}\bm{\gamma}(\bm{X})]\big)^{-1} \bm{H}.
$

In general, it is not possible to derive a simple natural gradient update formula for \eqref{eq:shift-form2} like \eqref{eq:updateUnconstQP}, \eqref{eq:update.formula.poly.on.rect}, and \eqref{eq:update.formula.poly.on.unit}, since the number of parameters is now less than that of random variables, and Theorem \ref{prop:char.exp.family_natgrad} cannot apply. However, the EM scheme is still applicable by minimizing a nonlinear convex optimization \eqref{eq:EM-step} at each iteration as far as an exponential family distribution with $A(\bm{\theta})-\mathbb{E}_{q_{\bm{x}}(\bm{\theta}^{(t-1)})}[\bm{\gamma}(\bm{X})]^\top\bm{\eta}(\bm{\theta})$ being convex is employed. 
Although the subproblem \eqref{eq:EM-step} might be as hard as (or harder than) the original problem, we may 
apply a generalized EM algorithm (GEM \cite{dempster1977maximum}), in which the M-step is terminated soon after a better incumbent is obtained. A straightforward example of such kind is the EM gradient algorithm \cite{lange1995gradient}, where a Newton method is applied to the subproblem. See Appendix \ref{sec:GEM} for the detailed descriptions.

\subsection{Minimization of a polynomial over a polytope}
\label{sec:poly.on.truncpoly}
The initial example of the extension \eqref{eq:shift-form2} is the generalized form of the problems \eqref{box_multi_problem} and \eqref{min_var_problem2}:
\begin{align}
\begin{array}{|ll}
\underset{\bm{\nu}\in\mathbb{R}^{p}}{\mbox{minimize}} \quad & F(\bm{\nu}):=\sum\limits_{i\in[I]}a_i \prod\limits_{j\in[p]}\nu_{j}^{n_{i,j}} \\
\mbox{subject to} & \bm{B}\bm{\nu}=\bm{c},~\bm{l}\leq\bm{\nu}\leq\bm{u},
\end{array}
\label{poly.on.truncpoly}
\end{align}
where $\bm{B} \in \mathbb{R}^{(p-n) \times p}$, $\bm{c} \in \mathbb{R}^{p-n}$ for $n,p\in\mathbb{N}$ such that $p>n$, and $\bm{l},\bm{u}\in\mathbb{R}^{p}$ such that $\bm{u}-\bm{l}\in\mathbb{R}^p_{++}$. 
We assume that \eqref{poly.on.truncpoly} has an interior feasible point. 
For example, the dual form of Support Vector Classification can be of the form \eqref{poly.on.truncpoly}. 
We first show that if $\bm{B}$ has full row rank, i.e., $\mathrm{rank}\bm{B}=p-n$, \eqref{poly.on.truncpoly} can be equivalently rewritten in the form of \eqref{eq:shift-form2} with $\Theta=[\bm{0},\bm{1}]$ so that the EM scheme based on independent binomial distributions can be developed. 
Similarly to Section \ref{sec:poly.on.rect}, the change of variables $\lambda_j=(\nu_j-l_j)/(u_j-l_j)$ for $j\in[p]$ 
yields an equivalent formulation to \eqref{poly.on.truncpoly} as
\begin{align}
\begin{array}{|ll}
\underset{\bm{\lambda}\in\mathbb{R}^{p}}{\mbox{minimize}} \quad & 
\phi(\bm{\lambda})=-\ln \Big(K-F\big(\mathrm{diag}(\bm{u}-\bm{l})\bm{\lambda}+\bm{l}\big)\Big) \\
\mbox{subject to} & \bm{B}\mathrm{diag}(\bm{u}-\bm{l})\bm{\lambda}=\bm{c}-\bm{B}\bm{l},~\bm{0}\leq\bm{\lambda}\leq\bm{1}
\end{array}
\label{poly.on.truncpoly_equiv}
\end{align}
for an upper bound $K$ of $F(\mathrm{diag}(\bm{u}-\bm{l})\bm{\lambda}+\bm{l})$ over the feasible region of \eqref{poly.on.truncpoly_equiv}. 
Without loss of generality, denoting $\bm{B}=\left( \bm{B}_1,\bm{B}_2 \right)$ for $\bm{B}_1 \in \mathbb{R}^{(p-n) \times n},\bm{B}_2 \in \mathbb{R}^{(p-n) \times (p-n)}$, we assume that $\bm{B}_2$ is nonsingular. 
The equality constraints 
of \eqref{poly.on.truncpoly_equiv} can then be rewritten as 
$\bm{\lambda}=\bm{H}\bm{\theta}+\bm{w},$ 
where 
\begin{align*}
\bm{H}&:= \mathrm{diag}(\bm{u}-\bm{l})^{-1}
\begin{pmatrix}
 \bm{I} \\
 - \bm{B}_2^{-1} \bm{B}_1
\end{pmatrix} \in \mathbb{R}^{p \times n},~
\bm{w}:=\mathrm{diag}(\bm{u}-\bm{l})^{-1}\left(
\begin{pmatrix}
 \bm{0} \\
\bm{B}_2^{-1} \bm{c}
\end{pmatrix}
-\bm{l}
\right)
 \in \mathbb{R}^p,
\end{align*}
and \eqref{poly.on.truncpoly} can be  
rewritten in the form of \eqref{eq:shift-form2} with $\Theta=[\bm{0},\bm{1}]$, as
\begin{align}
\begin{array}{|ll}
\underset{\bm{\lambda}\in\mathbb{R}^{p},\bm{\theta}\in\mathbb{R}^{n}}{\mbox{minimize}} \quad & \phi(\bm{\lambda})=
-\ln\Big(K-\big(\sum\limits_{i\in[I]}\tilde{a}_i \prod\limits_{j\in[p]}\lambda_{j}^{n_{i,j}}+\tilde{b}\big)\Big) \\
\mbox{subject to} & \bm{\lambda}=\bm{H}\bm{\theta}+\bm{w},~\bm{0}\leq\bm{\lambda}\leq\bm{1},
\end{array}
\label{poly.on.truncpoly2}
\end{align}
for some 
$\tilde{a}_i,~i\in [I]$, and $\tilde{b}$ 
such that  
\begin{align*}
 F(\mathrm{diag}(\bm{u}-\bm{l})\bm{\lambda}+\bm{l}) & = \sum_{i\in[I]} a_i \prod_{j\in[p]} \big((u_j-l_j)\lambda_j+l_j\big)^{n_{i,j}} 
 =  \sum_{i\in[I]} \tilde{a}_i \prod_{j\in[p]} \lambda_j^{n_{i,j}} + \tilde{b}.
\end{align*}

To apply the EM scheme to \eqref{poly.on.truncpoly2}, 
let $\bm{X}=(X_1,...,X_p)$ be a vector of $p$ independent binomial random variables, $X_j\sim\mathrm{Bin}(m_j,\lambda_j)$, with $m_j \geq \max\limits_{i\in[I]}
n_{i,j},~j \in [p]$, so that 
\begin{align}
p_{\bm{x}}(\bm{\lambda})
 &=\prod_{j\in[p]}\frac{m_j!}{x_j! (m_j-x_j)!}\lambda_j^{x_j}(1-\lambda_j)^{m_j-x_j}
~\big(=\hat{p}_{\bm{x}}(\bm{\theta}):=p_{\bm{x}}(\bm{H}\bm{\theta}+\bm{w})\big).
\label{eq:binomial.prob_shift}
\end{align}
Similarly to Proposition \ref{prop:expectation_condition_for_simple_binom} in Section \ref{sec:poly.on.rect}, the following proposition shows a function $G(\cdot)$ that satisfies the modified condition \eqref{eq:G-condition_ext} for problem \eqref{poly.on.truncpoly2} under the binomial probability distribution \eqref{eq:binomial.prob_shift}.
\begin{proposition}
Let $g_{i,j}(x_j)$ be defined as in Proposition \ref{prop:expectation_condition_for_simple_binom}. 
If $G$ is defined by
\begin{align}
G(\bm{x}) & = K - \Big( \sum_{i\in[I]} \tilde{a}_i \prod_{j\in[p]} g_{i,j}(x_j) + \tilde{b} \Big),
\tag{\ref{def:G}}
\end{align}
the condition \eqref{eq:G-condition_ext} for \eqref{poly.on.truncpoly2} is fulfilled, i.e., 
$f(\bm{\theta})(=\phi(\bm{H}\bm{\theta}+\bm{w}))=-\ln\mathbb{E}_{\mathbb{P}}\big[G(\bm{X})\big\vert\bm{\theta}\big]$ holds with $\hat{p}_{\bm{x}}(\bm{\theta})$ defined by \eqref{eq:binomial.prob_shift} and $\Theta=[\bm{0},\bm{1}]
$, or 
$\phi(\bm{\lambda})=-\ln\mathbb{E}_{\mathbb{P}}\big[G(\bm{X})\big\vert\bm{\lambda}\big]$ holds with $
p_{\bm{x}}(\bm{\lambda})$ defined by \eqref{eq:binomial.prob_shift} 
and $\Theta=[\bm{0},\bm{1}]$. 
\end{proposition}
For $\bm{\theta},\bar{\bm{\theta}}\in\mathrm{int}\{\bm{\theta}\in\mathbb{R}^n| \bm{H}\bm{\theta}+\bm{w}\in\Theta\}$, the surrogate function 
$u(\bm{\theta}\vert\bar{\bm{\theta}})$
 is a convex function of the form:
\begin{align*}
u(\bm{\theta}\vert\bar{\bm{\theta}})
&= - \sum\limits_{j\in[p]} 
\left\{
\mathbb{E}_{\hat{q}_{\bm{x}}(\bar{\bm{\theta}})}[X_j] \ln\frac{\bm{h}_j\bm{\theta}+w_j}{\bm{h}_j\bar{\bm{\theta}}+w_j}
+
(m_j-\mathbb{E}_{\hat{q}_{\bm{x}}(\bar{\bm{\theta}})}[X_j]) \ln\frac{1-(\bm{h}_j\bm{\theta}+w_j)}{1-(\bm{h}_j\bar{\bm{\theta}}+w_j)}
\right\}
+ f(\bar{\bm{\theta}}),
\end{align*}
where $\bm{h}_j$ denotes $j$-th row of $\bm{H}$, $\mathbb{E}_{\hat{q}_{\bm{x}}(\bar{\bm{\theta}})}[\bm{X}] = \mathrm{diag}(\bar{\bm{\lambda}})
 \big(\bm{m}  - 
 \mathrm{diag}(\bm{1}-\bar{\bm{\lambda}})\,
 \nabla\phi(\bar{\bm{\lambda}})\big) .$ 
with $\bar{\bm{\lambda}}=\bm{H}\bar{\bm{\theta}}+\bm{w}$ via \eqref{eq:Expectation_q}, and the EM algorithm for \eqref{poly.on.truncpoly2} is thus described as Algorithm \ref{alg:EMforPoly.Polytope:0}. 
In general, however, simple update formula is no longer available unlike the examples in Section \ref{sec:examples}.  
\begin{algorithm}[tb]
\caption{EM algorithm for \eqref{poly.on.truncpoly2}}
\label{alg:EMforPoly.Polytope:0}
\begin{algorithmic}[1]
\State 
Let $\bm{\theta}^{(0)}\in\mathbb{R}^n$ satisfy 
$\bm{H}\bm{\theta}^{(0)}+\bm{w}\in(\bm{0},\bm{1})$, 
and set $t\gets 1$.
\Repeat
\State $\bm{\lambda}^{(t-1)}\gets\bm{H}\bm{\theta}^{(t-1)}+\bm{w}$
\State $\mathbb{E}_{\hat{q}_{\bm{x}}(\bm{\theta}^{(t-1)})}[\bm{X}]\gets\mathrm{diag}(\bm{\lambda}^{(t-1)})\Big(\bm{m} -\mathrm{diag}\big(\bm{1}-\bm{\lambda}^{(t-1)}\big) \nabla\phi(\bm{\lambda}^{(t-1)})\Big)$,
\State
$
\bm{\theta}^{(t)}\gets
\arg\min\limits_{\bm{\theta}\in\mathbb{R}^{n}}\Big\{
- \sum\limits_{j\in[p]} \Big[\mathbb{E}_{\hat{q}_{\bm{x}}(\bm{\theta}^{(t-1)})}[X_j] \ln\frac{\bm{h}_j\bm{\theta}+w_j}{1-(\bm{h}_j\bm{\theta}+w_j)}
+m_j \ln(1-(\bm{h}_j\bm{\theta}+w_j))
\Big]
\Big\},$
\State and let $t\gets t+1$.
\Until{a termination condition is fulfilled}
\end{algorithmic}
\end{algorithm}
With $\bm{R}_1:= \mathrm{diag}(r_1,\ldots,r_n),\bm{R}_2:= \mathrm{diag}(r_{n+1},\ldots,r_p)$ where $r_j := \Big( \frac{\mathbb{E}_{\hat{q}_{\bm{x}}(\bar{\bm{\theta}})}[X_j]}{\bar{\lambda}_j^2} 
 + \frac{m_j-\mathbb{E}_{\hat{q}_{\bm{x}}(\bar{\bm{\theta}})}[X_j]}{(1-\bar{\lambda}_j)^2} \Big)^{-1}$,  $j \in[p]$,  
the Hessian of $u(\bm{\theta} | \bar{\bm{\theta}})$ with respect to $\bm{\theta}$ is given by 
\begin{align*}
\left.\nabla^2 u(\bm{\theta} | \bar{\bm{\theta}})\right|_{\bm{\theta}=\bar{\bm{\theta}}} 
& 
 = \bm{R}_1^{-1} + \left(\bm{B}_2^{-1}\bm{B}_1\right)^\top \bm{R}_2^{-1} \bm{B}_2^{-1}\bm{B}_1,
\end{align*}
and the Sherman--Morrison--Woodbury formula yields the inverse of the Hessian as 
\begin{align*}
\left(\nabla^2 u(\bm{\theta} | \bar{\bm{\theta}})\big\vert_{\bm{\theta}=\bar{\bm{\theta}}}\right)^{-1}
& = \bm{R}_1 - \bm{R}_1 \left(\bm{B}_2^{-1}\bm{B}_1\right)^\top \Big(\bm{R}_2 + \left(\bm{B}_2^{-1}\bm{B}_1\right) \bm{R}_1 \left(\bm{B}_2^{-1}\bm{B}_1\right)^\top \Big)^{-1} \bm{B}_2^{-1}\bm{B}_1 \bm{R}_1.
\end{align*}
This can be beneficial for efficient implementation of the EM gradient when the size of $\bm{B}_2$, i.e., $p-n$ is small.


\subsection{Quadratic program and linear program}
\label{subsec:QP}
The second example of the extension \eqref{eq:shift-form2} is 
a quadratic program (QP):
\begin{align}
\begin{array}{|ll}
\underset{\bm{\theta}\in\mathbb{R}^{n}}{\mbox{minimize}} \quad & \frac{1}{2}\bm{\theta}^\top\bm{Q}\bm{\theta}+\bm{b}^\top\bm{\theta} \\
\mbox{subject to} & \bm{A}^\top\bm{\theta}\leq\bm{c}, \\
\end{array}
\label{eq:quadratic_problem}
\end{align}
where $\bm{A}\in\mathbb{R}^{n \times p},\bm{b} \in \mathbb{R}^n,\bm{c}\in\mathbb{R}^p$, and $\bm{Q}\in\mathbb{S}^{n}$, so the objective can be nonconvex. It is easy to see that \eqref{eq:quadratic_problem} falls in the form \eqref{eq:shift-form2} with $\bm{\lambda}=\bm{c}-\bm{A}^\top\bm{\theta}$ and $\Theta=\mathbb{R}^p_+$.  
Moreover, \eqref{eq:quadratic_problem} includes a linear program (LP) as a special case by setting $\bm{Q}=\bm{O}$. We assume that \eqref{eq:quadratic_problem} has an interior point. 

To apply the EM algorithm to 
\eqref{eq:quadratic_problem}, consider $p+n$ random variables $\bm{X}:=(\bm{X}_1,\bm{X}_2)$ such that $\bm{X}_1=(X_{1,1},...,X_{1,p})$ is $p$ independent Poisson random variables, $\bm{X}_2=(X_{2,1},...,X_{2,n})$ follows $n$-dimensional normal distribution $\bm{\mathrm{N}}(\bm{\theta},\bm{\Sigma})$ for $\bm{\Sigma}\in\mathbb{S}^{n}_{++}$, and every pair of $X_{1,i}$ and $X_{2,j}$ is independent for $i \in [p],j\in [n]$. Specifically, its likelihood-like function of $\bm{X}=\bm{x}$ is given by
\begin{align}
p_{\bm{x}}(\bm{\theta}) & = p_{\bm{x}_1}(\bm{\lambda})\,p_{\bm{x}_2}(\bm{\theta})
= p_{\bm{x}_1}(\bm{c}-\bm{A}^\top\bm{\theta})\,p_{\bm{x}_2}(\bm{\theta})
= \hat{p}_{\bm{x}_1}(\bm{\theta})\,p_{\bm{x}_2}(\bm{\theta})
,
\label{eq:Poisson+Normal0}
\end{align}
where 
\begin{align*}
p_{\bm{x}_1}(\bm{\lambda}) 
&:=\prod_{j\in[p]}\frac{\lambda_{j}^{x_{1,j}}\exp(-\lambda_{j})}{x_{1,j}!} =\underbrace{\Big(\prod\limits_{j\in[p]}x_{1,j}!\Big)^{-1}}_{=\beta(\bm{x}_1)}
\exp\Big(
\underbrace{\sum_{j\in[p]}x_{1,j}\ln\lambda_j}_{=\bm{\gamma}(\bm{x}_1)^\top\bm{\eta}(\bm{\lambda})}
\underbrace{-
\sum_{j\in[p]}\lambda_j}_{=-A(\bm{\lambda})}\Big)
\\
&=\hat{p}_{\bm{x}_1}(\bm{\theta}) 
:=\prod_{j\in[p]}\frac{\big(c_j - \bm{a}_j^\top \bm{\theta}\big)^{x_{1j}}\exp\big(-(c_j - \bm{a}_j^\top \bm{\theta})\big)}{x_{1,j}!}, \\
p_{\bm{x}_2}(\bm{\theta}) 
& :=\frac{1}{(2\pi)^{\frac{n}{2}} \sqrt{|\bm{\Sigma}|}} \exp\Big(- \frac{1}{2}(\bm{x}_2-\bm{\theta})^\top\bm{\Sigma}^{-1} (\bm{x}_2-\bm{\theta}) \Big),
\end{align*}
where $\bm{a}_j$ and $\bm{q}_j$ are the $j$-th column of $\bm{A}$ and $\bm{Q}$, respectively, 
i.e., $\bm{A}=(\bm{a}_1,\ldots,\bm{a}_p)$ and $\bm{Q}=(\bm{q}_1,\ldots,\bm{q}_{n})$. 
It is easy to see that \eqref{eq:Poisson+Normal0} is a member of the exponential family.
\begin{proposition}
\label{prop:G.for.quadratic_problem0}
If $\bm{\Sigma}^{-1}\succ\bm{Q}$ holds and $G$ is defined by $G(\bm{x}_1,\bm{x}_2)= G_1(\bm{x}_1)G_2(\bm{x}_2),$ 
where
\begin{align*}
&G_1(\bm{x}_1) = 1, \quad G_2(\bm{x}_2) 
= \sqrt{\frac{|\bm{\Sigma}|}{|\tilde{\bm{Q}} + \bm{\Sigma}^{-1}|}} \exp \Big( - \frac{1}{2}\bm{x}_2^\top \tilde{\bm{Q}}\bm{x}_2 - \tilde{\bm{b}}^\top \bm{x}_2 \Big),\\
&\tilde{\bm{Q}} := \bm{\Sigma}^{-1}(\bm{\Sigma}^{-1}-\bm{Q})^{-1} \bm{\Sigma}^{-1} - \bm{\Sigma}^{-1},
 \quad\mbox{and}\quad
 \tilde{\bm{b}} 
 := \bm{\Sigma}^{-1}(\bm{\Sigma}^{-1}-\bm{Q})^{-1} \bm{b},
\end{align*}
the condition \eqref{eq:G-condition_ext} for \eqref{eq:quadratic_problem} is fulfilled 
with $
p_{\bm{x}}(\bm{\theta})$ defined by \eqref{eq:Poisson+Normal0}, $f(\bm{\theta}) = \frac{1}{2}\bm{\theta}^\top\bm{Q}\bm{\theta}+\bm{b}^\top\bm{\theta}$, and $\Theta=\mathbb{R}^p_+$.
\end{proposition}
See Appendix \ref{sec:proof:prop:EMcond_for_QP0} for the proof. 

For $\bm{\theta},\bar{\bm{\theta}}\in
\mathrm{int}\{\bm{\theta}\vert\bm{A}^\top\bm{\theta}\leq\bm{c}\}$, 
the surrogate $u(\bm{\theta} | \bar{\bm{\theta}})$ is a convex function of the form:
\begin{align*}
u(\bm{\theta} | \bar{\bm{\theta}}) & = -\sum\limits_{j\in[p]} \Big\{(c_j-\bm{a}_j^\top\bar{\bm{\theta}}) 
\ln\frac{c_j-\bm{a}_j^\top\bm{\theta}}{c_j-\bm{a}_j^\top\bar{\bm{\theta}}}+\bm{a}_j^\top(\bm{\theta}-\bar{\bm{\theta}})\Big\} \\
& \qquad - \big\{ 
(\bm{\Sigma}^{-1}-\bm{Q})\bar{\bm{\theta}}-\bm{b}
\big\}^\top
\big(\bm{\theta}-\bar{\bm{\theta}}\big) + \frac{1}{2}\bm{\theta}^\top\bm{\Sigma}^{-1}\bm{\theta} - \frac{1}{2}\bar{\bm{\theta}}^\top\bm{\Sigma}^{-1}\bar{\bm{\theta}} + f(\bar{\bm{\theta}}),
\end{align*}
and the EM algorithm for \eqref{eq:quadratic_problem} is described in Algorithm \ref{alg:EMforQP0}.
\begin{algorithm}[tb]
\caption{EM algorithm for quadratic program \eqref{eq:quadratic_problem}}
\label{alg:EMforQP0}
\begin{algorithmic}[1]
\State 
Let $\bm{\theta}^{(0)}\in\mathbb{R}^{n}$ be an interior point of the feasible region  of \eqref{eq:quadratic_problem}, and set $t\gets 1$. 
\Repeat
\State
 For $\mu^{(t-1)}_j:= c_j-\bm{a}_j^\top\bm{\theta}^{(t-1)}$ and $\bm{\nu}^{(t-1)}:=\bm{b}-(\bm{\Sigma}^{-1}-\bm{Q})\bm{\theta}^{(t-1)}$, solve
\begin{align}
\bm{\theta}^{(t)}\gets
\arg\min_{\bm{\theta}\in\mathbb{R}^{p}
}
\Big\{
\begin{array}{r}
-\sum\limits_{j\in[p]} \Big( \bm{a}_j^\top\bm{\theta}  + \mu^{(t-1)}_j \ln(c_j-\bm{a}_j^\top\bm{\theta}) \Big)
+ \frac{1}{2}\bm{\theta}^\top\bm{\Sigma}^{-1}\bm{\theta}+(\bm{\nu}^{(t-1)})^\top\bm{\theta}
\end{array}
\Big\},
\label{eq:M-step_for_QP}
\end{align}
\State and let $t\gets t+1$.
\Until{a termination condition is fulfilled}
\end{algorithmic}
\end{algorithm}
When the objective function is convex, the convergence of Algorithm \ref{alg:EMforQP0} can be shown based on \cite{teboulle1997convergence} (see Section \ref{proof:prop:QPconvergence} for the outline of its proof). 
\begin{proposition}\label{prop:QPconvergence}
Suppose that $\bm{Q}\in\mathbb{S}^{n}_+$, $\bm{A}$ has full row rank, i.e., $\mathrm{rank}A=n$, and \eqref{eq:quadratic_problem} has an interior point and optimum. Then the sequence, $\{\bm{\theta}^{(t)}\}_{t\geq 0}$, generated by Algorithm \ref{alg:EMforQP0} converges to an optimal solution $\bm{\theta}^{\star}$ of \eqref{eq:quadratic_problem}. Moreover, if $X^\star:=\{\bm{\theta}\vert f(\bm{\theta})=f(\bm{\theta}^\star)\}\neq\emptyset$, the sequence $\{\bm{\theta}^{(t)}\}_{t\geq 0}$ converges to an optimal solution of \eqref{eq:quadratic_problem}. 
\end{proposition}

While the MM algorithm can be approached by a GEM, it may be possible to derive simple algorithms if the polyhedral feasible set is given by a simple one. For example, consider the box constrained (possibly, nonconvex) quadratic minimization:
\begin{align}
\begin{array}{|ll}
\underset{\bm{\theta}\in\mathbb{R}^n}{\mbox{minimize}} & \frac{1}{2}\bm{\theta}^\top\bm{Q}\bm{\theta} + \bm{b}^\top\bm{\theta} \\
\mbox{subject to} & \bm{l}\leq\bm{\theta}\leq\bm{u}, \\
\end{array}
\label{quadratic_problem_box2}
\end{align}
which is of the form \eqref{eq:quadratic_problem} with $\bm{c}=(\bm{u}^\top,-\bm{l}^\top)^\top$ and $\bm{A}=(\bm{I}, ~-\bm{I})$. 
For $\bm{\theta}^{(t-1)}=\bar{\bm{\theta}}$, the $\bm{\theta}^{(t)}$ is obtained as the solution $\bm{\theta}$ of
\begin{align*}
\nabla u(\bm{\theta} | \bar{\bm{\theta}}) 
& = - \big[\mathrm{diag}(\bm{\theta}-\bm{l})\big]^{-1}(\bar{\bm{\theta}}-\bm{l}) + 
\big[\mathrm{diag}(\bm{u}-\bm{\theta})\big]^{-1}(\bm{u}-\bar{\bm{\theta}}) - (\bm{\Sigma}^{-1}-\bm{Q})\bar{\bm{\theta}} +\bm{b} + \bm{\Sigma}^{-1}\bm{\theta}=\bm{0}. 
\end{align*}
If $\bm{\Sigma}$ is a diagonal matrix $\mathrm{diag}(\bm{\sigma})$ for $\sigma_j>0,~j\in[n]$, 
the $j$-th component of $\bm{\theta}^{(t)}$ is obtained 
as a solution satisfying $\theta_j\in(l_j,u_j)$ and 
\begin{align}
-\frac{\bar{\theta}_j-l_j}{\theta_j-l_j}+\frac{u_j-\bar{\theta}_j}{u_j-\theta_j}-\kappa_j+\sigma_j^{-1}\theta_j=0,\label{box_eq_const}
\end{align}
where $\bar{\theta}_j=\theta^{(t-1)}_j$ and $\kappa_j:=[(\bm{\Sigma}^{-1}-\bm{Q})\bar{\bm{\theta}}-\bm{b}]_j,~j\in[n]$. 
Note that \eqref{box_eq_const} has a unique solution since the left-hand side increases monotonically from $-\infty$ to $\infty$ on the interval $(l_j,u_j)$, and it is easy to see that the solution is obtained via the explicit solution formula for the cubic equation:
\begin{align*}
\theta_j^3-(\sigma_j\kappa_j+l_j+u_j)\theta_j^2+[l_ju_j+(l_j+u_j)(\kappa_j-1)\sigma_j]\theta_j+\sigma_j[\bar{\theta}_j(u_j-l_j)-\kappa_jl_ju_j]=0,
\end{align*}
which is obtained by multiplying \eqref{box_eq_const} by $-\sigma_j(\theta_j-l_j)(u_j-\theta_j)$.
Accordingly, the MM algorithm results in a simple component-wise formula like \eqref{eq:update.formula.poly.on.rect} and \eqref{eq:update.formula.poly.on.unit}. 
\eqref{quadratic_problem_box2} can also be approached using \eqref{eq:update.formula.poly.on.rect}, but it is interesting to note that this one is different from \eqref{eq:update.formula.poly.on.rect}. 
\begin{remark}
We can construct another EM/MM algorithm based on the combination of the Poisson and normal distributions for another form of QP, which is explained in Appendix \ref{subsec:QP2} since the construction is parallel to the one for \eqref{eq:quadratic_problem}. On the other hand, using this derivation, the EM/MM algorithm for LP can be constructed based only on the Poisson and does not require the normal distribution. See the appendix for the details. 
\end{remark}


\section{Concluding remarks}\label{sec:conclusion}
In this paper we develop an algorithmic scheme to approach general optimization problems on the basis of the EM algorithm. 
We demonstrate that when a member of the exponential family distributions is used for some classes of optimization problems with polynomial objectives, it falls into a natural gradient descent algorithm of the employed probability distribution with a constant step size, or an iterative nonlinear convex optimization. 

The application of EM algorithm to general optimization problems opens up a door to developing a new class of approaches to optimization problems, but developing efficient algorithms and expanding its utility beyond their aesthetic appeal are left for future researches. 

In particular, while we collectively present several linkages between optimization problems and probability distributions, some problems can be approached via different reformulations and distributions. For example, 
the quadratic minimization over the rectangle can be solved not only via the approach (ii) in Table \ref{tbl:problems} but also (iv) and (vi), as demonstrated in the final paragraph in the previous section. Further comparison among different formulations or that among different distributions are also to be investigated. 

Relations to information geometry may be an interesting research topic as well. The
relation to interior point algorithms has been occasionally pointed out in the literature  \cite{tanabe1980geometric,ohara1999information,kakihara2013information}, but the analysis based on the exponential family distributions is not examined yet except in the context of statisitical EM algorithms (e.g., \cite{amari2001differential,cousseau2008dynamics}). 
In information geometry, the EM algorithm is understood as optimization on a manifold based on data (e.g., \cite{amari1995information}). 
%
While a geometric interpretation of the EM algorithm for some estimation problems with probabilistic models has been studied (e.g., \cite{6796330} for learning hidden variable models in neural networks), 
it remains to be seen whether these existing perspectives, in conjunction with the EM algorithm developed in this paper, can help to discover new insights into general (non-statistical) optimization problems.

\paragraph{Acknowledgment}
The second author's research is supported in part by MEXT Grant-in-Aid 
24K01113.
The careful reading and constructive comments of two reviewers have greatly contributed to improving the readability of this paper. Needless to say, any unintentional errors remaining in the final draft are the responsibility of the authors. 

\bibliography{references}
\bibliographystyle{plain}

\begin{appendices}

\section{Derivation of formulae and proof of propositions}

\subsection{Proof of Proposition \ref{prop:char.exp.family}.}\label{sec:proof.prop:char.exp.family}
\noindent{\it (a):} The claim is confirmed by simply substituting $\tilde{A}(\bm{\theta}^{(t-1)})$ and $\tilde{\beta}(\bm{x})$ into the formula to be proved.
\\\noindent{\it (b):}  
By the definition \eqref{eq:def_u} of the surrogate function $u$, we have 
\begin{align*}
\lefteqn{u(\bm{\theta}\vert\bm{\theta}^{(t-1)})}\\
 & = - \mathbb{E}_{q_{\bm{x}}(\bm{\theta}^{(t-1)})}[\ln G(\bm{X})] - \mathbb{E}_{q_{\bm{x}}(\bm{\theta}^{(t-1)})}[\ln p_{\bm{X}}(\bm{\theta})] + \mathbb{E}_{q_{\bm{x}}(\bm{\theta}^{(t-1)})}[\ln q_{\bm{X}}(\bm{\theta}^{(t-1)})] \\
 & = - \mathbb{E}_{q_{\bm{x}}(\bm{\theta}^{(t-1)})}[\ln G(\bm{X})] - \mathbb{E}_{q_{\bm{x}}(\bm{\theta}^{(t-1)})}[\ln p_{\bm{X}}(\bm{\theta})] \\
 & \quad + \big( \mathbb{E}_{q_{\bm{x}}(\bm{\theta}^{(t-1)})}[\ln G(\bm{X})] + \mathbb{E}_{q_{\bm{x}}(\bm{\theta}^{(t-1)})}[\ln p_{\bm{X}}(\bm{\theta}^{(t-1)})] + f(\bm{\theta}^{(t-1)}) \big) & (\because \eqref{eq:G-condition},\eqref{eq:q_x}) \\
 & = - \mathbb{E}_{q_{\bm{x}}(\bm{\theta}^{(t-1)})}[\ln p_{\bm{X}}(\bm{\theta})] + \mathbb{E}_{q_{\bm{x}}(\bm{\theta}^{(t-1)})}[\ln p_{\bm{X}}(\bm{\theta}^{(t-1)})] + f(\bm{\theta}^{(t-1)}) \\
& =  - \mathbb{E}_{q_{\bm{x}}(\bm{\theta}^{(t-1)})}[\bm{\gamma}(\bm{X})]^\top\big(\bm{\eta}(\bm{\theta})-\bm{\eta}(\bm{\theta}^{(t-1)})\big) + A(\bm{\theta}) - A(\bm{\theta}^{(t-1)}) + f(\bm{\theta}^{(t-1)}),& \tag{\ref{eq:u_component}}
\end{align*}
which shows the formula \eqref{eq:u_component}.
\\\noindent{\it (c):} 
Although the statement (c) is known in the context of EM algorithms, we here present its proof to make the article self-contained. 
Note 
\begin{align}
 \nabla A(\bm{\theta})
 & = \frac{1}{\exp(A(\bm{\theta}))} \nabla\Big(\int_{\bm{x}\in\Xi}\beta(\bm{x})\exp\big(\bm{\eta}(\bm{\theta})^\top\bm{\gamma}(\bm{x})\big) \mathrm{d}\bm{x}\Big) & (\because \eqref{eq:expdist})\nonumber \\
 & = \nabla\bm{\eta}(\bm{\theta})^\top\int_{\bm{x}\in\Xi}\bm{\gamma}(\bm{x})\beta(\bm{x})\exp\big(\bm{\eta}(\bm{\theta})^\top\bm{\gamma}(\bm{x})-A(\bm{\theta})\big) \mathrm{d}\bm{x}\nonumber \\
 & = \nabla\bm{\eta}(\bm{\theta})^\top \mathbb{E}_{p_{\bm{x}}(\bm{\theta})}[\bm{\gamma}(\bm{X})],
 \label{eq:nabla.ln.alpha}
\end{align}
where all the gradient operations $\nabla$ apply with respect to $\bm{\theta}$. 
Owing to \eqref{eq:u_component} and \eqref{eq:nabla.ln.alpha}, we have
\begin{align*}
 \nabla_{\bm{\theta}}u(\bm{\theta}|\bm{\theta}^{(t-1)})&=\nabla A(\bm{\theta})-\nabla\bm{\eta}(\bm{\theta})^\top\mathbb{E}_{q_{\bm{x}}(\bm{\theta}^{(t-1)})}[\bm{\gamma}(\bm{X})] & (\because \eqref{eq:u_component})\\
  & = \nabla\bm{\eta}(\bm{\theta})^\top \big( \mathbb{E}_{p_{\bm{x}}(\bm{\theta})}[\bm{\gamma}(\bm{X})] - \mathbb{E}_{q_{\bm{x}}(\bm{\theta}^{(t-1)})}[\bm{\gamma} (\bm{X})] \big), & (\because \eqref{eq:nabla.ln.alpha})
\end{align*}
which leads immediately to the claim.
\\\noindent{\it (d):} 
Observe first that \eqref{eq:nabla.ln.alpha} holds true also for $\tilde{A}(\bm{\theta})$ and $q_{\bm{x}}(\bm{\theta})$ in place of $A(\bm{\theta})$ and $p_{\bm{x}}(\bm{\theta})$, respectively, i.e.,
\begin{align}
\nabla\tilde{A}(\bm{\theta}) = \nabla\bm{\eta}(\bm{\theta})^\top \mathbb{E}_{q_{\bm{x}}(\bm{\theta})}[\bm{\gamma}(\bm{X})]. 
\label{eq:nabla.ln.alpha_hat}
\end{align}
Using this, we have
\begin{align*}
  \nabla f(\bm{\theta})
 & = \nabla \left(-\ln\mathbb{E}_{p_{\bm{x}}(\bm{\theta})}[G(\bm{X})]\right) \\
 & = \nabla \big(A(\bm{\theta})-\tilde{A}(\bm{\theta})\big) & (\because \eqref{q_exponential_family}) \\
& = - \nabla\bm{\eta}(\bm{\theta})^\top \mathbb{E}_{q_{\bm{x}}(\bm{\theta})}\big[ \bm{\gamma}(\bm{X}) \big] + \nabla\bm{\eta}(\bm{\theta})^\top \mathbb{E}_{p_{\bm{x}}(\bm{\theta})}[ \bm{\gamma}(\bm{X})], & (\because\eqref{eq:nabla.ln.alpha}\mbox{~and~}\eqref{eq:nabla.ln.alpha_hat})
\end{align*}
which leads to \eqref{eq:Expectation_q} by multiplying the inverse of $\nabla\bm{\eta}(\bm{\theta})^\top$.


\subsection{Proof of Corollary \ref{cor:FOS}}\label{sec:proof:cor:FOS}
The gradient of $h(\bm{\theta}\vert \bar{\bm{\theta}}) := u(\bm{\theta}\vert \bar{\bm{\theta}}) - f(\bm{\theta})$ in $\bm{\theta}$ is given as
\begin{align*}
 \nabla h(\bm{\theta} | \bar{\bm{\theta}})
 = \nabla\mathrm{KL}\left( q_{\bm{x}}(\bar{\bm{\theta}}) || q_{\bm{x}}(\bm{\theta}) \right) 
 & 
 = \nabla \Big( \int_{\bm{x}\in\Xi} q_{\bm{x}}(\bar{\bm{\theta}}) \ln{ \frac{ q_{\bm{x}}(\bar{\bm{\theta}}) }{ q_{\bm{x}}(\bm{\theta}) } } \mathrm{d}\bm{x} \Big) \\
 & = - \nabla \Big( \int_{\bm{x} \in\Xi } q_{\bm{x}}(\bar{\bm{\theta}}) \ln{ q_{\bm{x}}(\bm{\theta}) } \mathrm{d}\bm{x} \Big) \\
 & = - \nabla \Big( \int_{\bm{x} \in\Xi } q_{\bm{x}}(\bar{\bm{\theta}}) \bm{\eta}(\bm{\theta})^\top \bm{\gamma}(\bm{x}) \mathrm{d}\bm{x} \Big) + \nabla \Big( \int_{\bm{x} \in\Xi } q_{\bm{x}}(\bar{\bm{\theta}}) 
 \tilde{A}(\bm{\theta})
  \mathrm{d}\bm{x} \Big) \\
 & = - \int_{\bm{x} \in \Xi } \nabla \bm{\eta}(\bm{\theta})^\top \bm{\gamma}(\bm{x}) q_{\bm{x}}(\bar{\bm{\theta}}) \mathrm{d}\bm{x} + \nabla \tilde{A}(\bm{\theta}) \\
 & = \nabla \bm{\eta}(\bm{\theta})^\top \big( \mathbb{E}_{q_{\bm{x}}(\bm{\theta})}[ \bm{\gamma}(\bm{X})] - \mathbb{E}_{q_{\bm{x}}(\bar{\bm{\theta}})}[ \bm{\gamma}(\bm{X})]\big),
\end{align*}
which equals $0$ at $\bm{\theta}=\bar{\bm{\theta}}$.


\subsection{Proof of Theorem \ref{prop:char.exp.family_natgrad}}\label{sec:proof:prop:char.exp.family_natgrad}
Observe first that
\begin{align}
\nabla \ln p_{\bm{x}}(\bm{\theta}) & =  
\nabla \left( \ln \beta(\bm{x}) - A(\bm{\theta}) + \bm{\eta}(\bm{\theta})^\top \bm{\gamma}(\bm{x}) \right) & (\because \eqref{eq:expdist})\nonumber \\
 & = \nabla \big( \bm{\eta}(\bm{\theta})^\top \bm{\gamma}(\bm{x}) - A (\bm{\theta}) \big) & (\because \nabla\big(\ln \beta(\bm{x})\big)=\bm{0}) \nonumber \\
 & = \nabla\bm{\eta}(\bm{\theta})^\top \bm{\gamma}(\bm{x}) - \nabla\bm{\eta}(\bm{\theta})^\top \mathbb{E}_{p_{\bm{x}}(\bm{\theta})}[\bm{\gamma}(\bm{X})] & (\because\eqref{eq:nabla.ln.alpha}) \nonumber \\
 & = \nabla\bm{\eta}(\bm{\theta})^\top \bm{\Omega}^{-1} \left( \bm{\Omega} \bm{\gamma}(\bm{x}) - \mathbb{E}_{p_{\bm{x}}(\bm{\theta})}[\bm{\Omega} \bm{\gamma}(\bm{X})] \right) \nonumber \\
 & = \nabla\bm{\eta}(\bm{\theta})^\top \bm{\Omega}^{-1} \left( \bm{\Omega} \bm{\gamma}(\bm{x}) - \bm{\theta} \right).& (\because\eqref{eq:unbiased_estimator}) \label{eq:p_logL}
\end{align}
Note that 
\begin{align}
\mathbb{E}_{p_{\bm{x}}(\bm{\theta})}\Big[\big(\bm{\Omega}\bm{\gamma}(\bm{X}) - \bm{\theta}\big) \nabla \ln p_{\bm{X}}(\bm{\theta})^\top\Big]
& = \int_{\bm{x}\in\Xi} \Big\{ \big( \bm{\Omega} \bm{\gamma}(\bm{x}) - \bm{\theta} \big) \nabla \ln p_{\bm{x}}(\bm{\theta})^\top \Big\} p_{\bm{x}}(\bm{\theta})\,\mathrm{d}\bm{x} \nonumber \\
& = \int_{\bm{x}\in\Xi} \left( \bm{\Omega} \bm{\gamma}(\bm{x}) - \bm{\theta} \right) \nabla p_{\bm{x}}(\bm{\theta})^\top\,\mathrm{d}\bm{x} \nonumber \\
& = \bm{I} \label{eq:unbiased_likelihood}
\end{align}
since we have 
$
\int_{\bm{x}\in\Xi} \bm{\Omega} \bm{\gamma}(\bm{x}) \nabla p_{\bm{x}}(\bm{\theta})^\top\,\mathrm{d}\bm{x}=\bm{I}$ and $\int_{\bm{x}\in\Xi} \nabla p_{\bm{x}}(\bm{\theta})\,\mathrm{d}\bm{x}=\bm{0}
$ 
by differentiating \eqref{eq:unbiased_estimator} and $\int_{\bm{x}\in\Xi} p_{\bm{x}}(\bm{\theta})\,\mathrm{d}\bm{x}=1$ with respect to $\bm{\theta}$, respectively.  
Therefore, we have
\begin{align*}
\bm{I} & = \mathbb{E}_{p_{\bm{x}}(\bm{\theta})}[\left( \bm{\Omega} \bm{\gamma}(\bm{X}) - \bm{\theta} \right) \nabla \ln p_{\bm{X}}(\bm{\theta})^\top] & (\because\eqref{eq:unbiased_likelihood}) \\
& = \mathbb{E}_{p_{\bm{x}}(\bm{\theta})}[\left( \bm{\Omega} \bm{\gamma}(\bm{X}) - \bm{\theta} \right)\left( \bm{\Omega} \bm{\gamma}(\bm{X}) - \bm{\theta} \right)^\top ] \big( \nabla\bm{\eta}(\bm{\theta})^\top \bm{\Omega}^{-1} \big)^\top & (\because\eqref{eq:p_logL})  \\
& = \mathbb{V}_{p_{\bm{x}}(\bm{\theta})}[\bm{\Omega} \bm{\gamma}(\bm{X})] \big( \nabla\bm{\eta}(\bm{\theta})^\top \bm{\Omega}^{-1} \big)^\top, 
\end{align*}
so we have $\mathbb{V}_{p_{\bm{x}}(\bm{\theta})}[\bm{\Omega} \bm{\gamma}(\bm{X})] = \big(\nabla\bm{\eta}(\bm{\theta})^\top \bm{\Omega}^{-1} \big)^{-1}$, establishing the second equality of \eqref{fisher_inv}.
With this, the first equality of \eqref{fisher_inv} can also be confirmed as
\begin{align*}
 \bm{I}_{p_{\bm{x}}}(\bm{\theta})
 & = \nabla\bm{\eta}(\bm{\theta})^\top \bm{\Omega}^{-1} \mathbb{E}_{p_{\bm{x}}(\bm{\theta})} \big[ \left( \bm{\Omega} \bm{\gamma}(\bm{X}) - \bm{\theta} \right) \left( \bm{\Omega} \bm{\gamma}(\bm{X}) - \bm{\theta} \right)^\top \big] \nabla\bm{\eta}(\bm{\theta})^\top \bm{\Omega}^{-1} 
 & (\because \eqref{eq:p_logL}) \\
 & = \mathbb{V}_{p_{\bm{x}}(\bm{\theta})}[\bm{\Omega}\bm{\gamma}(\bm{X})]^{-1} 
 \mathbb{V}_{p_{\bm{x}}(\bm{\theta})}[\bm{\Omega}\bm{\gamma}(\bm{X})] 
 \, \mathbb{V}_{p_{\bm{x}}(\bm{\theta})}[\bm{\Omega}\bm{\gamma}(\bm{X})]^{-1}
  \\
 & 
 = \mathbb{V}_{p_{\bm{x}}(\bm{\theta})}[\bm{\Omega}\bm{\gamma}(\bm{X})]^{-1}.
\end{align*}
Recalling 
$
 \nabla_{\bm{\theta}}u(\bm{\theta}|\bm{\theta}^{(t-1)})
 = \nabla\bm{\eta}(\bm{\theta})^\top \big( \mathbb{E}_{p_{\bm{x}}(\bm{\theta})}[\bm{\gamma}(\bm{X})] - \mathbb{E}_{q_{\bm{x}}(\bm{\theta}^{(t-1)})}[\bm{\gamma} (\bm{X})] \big),
$ which appears in the proof of 
Proposition \ref{prop:char.exp.family}(c), we see that if $\bm{\theta}^{(t)}$ is a stationary point, i.e., $\nabla_{\bm{\theta}} u(\bm{\theta}^{(t)}\vert\bm{\theta}^{(t-1)})=\bm{0}$, then we have $\mathbb{E}_{p_{\bm{x}}(\bm{\theta}^{(t)})}[\bm{\gamma}(\bm{X})]=\mathbb{E}_{q_{\bm{x}}(\bm{\theta}^{(t-1)})}[\bm{\gamma} (\bm{X})]$ since $\nabla\bm{\eta}(\bm{\theta})$ is non-singular. 
For $\bm{\theta}=\bm{\theta}^{(t)}$, 
\eqref{eq:Expectation_q} becomes
\begin{align}
\mathbb{E}_{q_{\bm{x}}(\bm{\theta}^{(t)})}\big[\bm{\gamma}(\bm{X}) \big]
& = \mathbb{E}_{p_{\bm{x}}(\bm{\theta}^{(t)})}[\bm{\gamma}(\bm{X})] - \big(\nabla\bm{\eta}(\bm{\theta}^{(t)})^\top \big)^{-1} \nabla f(\bm{\theta}^{(t)}).
\label{eq:Eq=Ep-etaTgradf}
\end{align}
With this and \eqref{eq:unbiased_estimator}, we reach the formula \eqref{eq:natgrad} as follows.
\begin{align*}
\bm{\theta}^{(t)}
 & =\mathbb{E}_{p_{\bm{x}}(\bm{\theta}^{(t)})}[\bm{\Omega}\bm{\gamma}(\bm{X})] & (\because \eqref{eq:unbiased_estimator})\\
 & =\bm{\Omega} \mathbb{E}_{p_{\bm{x}}(\bm{\theta}^{(t)})}[\bm{\gamma}(\bm{X})] \\
 & = \bm{\Omega} \mathbb{E}_{q_{\bm{x}}(\bm{\theta}^{(t-1)})}[\bm{\gamma}(\bm{X})] & (\because \mbox{ Proposition }\ref{prop:char.exp.family}(c))\\
 & =\mathbb{E}_{p_{\bm{x}}(\bm{\theta}^{(t-1)})}[\bm{\Omega} \bm{\gamma}(\bm{X})] - \bm{\Omega} \big( \nabla\bm{\eta}(\bm{\theta}^{(t-1)})^\top \big)^{-1} \nabla f(\bm{\theta}^{(t-1)}) & (\because 
 \eqref{eq:Eq=Ep-etaTgradf} \mbox{ with }t\gets t-1 )\\
 & = \bm{\theta}^{(t-1)} - \big(\nabla\bm{\eta}(\bm{\theta}^{(t-1)})^\top \bm{\Omega}^{-1} \big)^{-1} \nabla f(\bm{\theta}^{(t-1)})&(\because \eqref{eq:unbiased_estimator} \mbox{ with }\bm{\theta}=\bm{\theta}^{(t-1)})\\
 & = \bm{\theta}^{(t-1)} - \bm{I}_{p_{\bm{x}}}(\bm{\theta}^{(t-1)})^{-1} \nabla f(\bm{\theta}^{(t-1)}). &(\because \eqref{fisher_inv})
\end{align*}


\subsection{Proof of Proposition \ref{prop:G.for.cvx.q.min}}\label{sec:proof:prop:G.for.cvx.q.min}
Since $\bm{\Sigma}^{-1}-\bm{Q}$ is positive definite, so is $\tilde{\bm{Q}}+\bm{\Sigma}^{-1}=\bm{\Sigma}^{-1}(\bm{\Sigma}^{-1}-\bm{Q})^{-1}\bm{\Sigma}^{-1}$, and we have $|\tilde{\bm{Q}}+\bm{\Sigma}^{-1}|>0$ and $(\tilde{\bm{Q}}+\bm{\Sigma}^{-1})^{-1}=\bm{\Sigma}(\bm{\Sigma}^{-1}-\bm{Q})\bm{\Sigma}$. 
We then see that
\begin{align*}
\mathbb{E}_{p_{\bm{x}}(\bm{\theta})}\big[G(\bm{X})\big]
& = \sqrt{\frac{|\bm{\Sigma}|}{|\tilde{\bm{Q}} + \bm{\Sigma}^{-1}|}} \int_{\bm{x} \in \mathbb{R}^p} \exp \left( - \frac{1}{2}\bm{x}^\top \tilde{\bm{Q}}\bm{x} - \tilde{\bm{b}}^\top \bm{x} \right) p_{\bm{x}}(\bm{\theta}) \mathrm{d}\bm{x} \\
& = \frac{\exp \left( - \frac{1}{2} \bm{\theta}^\top \bm{\Sigma}^{-1} \bm{\theta} \right)}{(2\pi)^{\frac{p}{2}} \sqrt{|\tilde{\bm{Q}} + \bm{\Sigma}^{-1}|}} \int_{\bm{x} \in \mathbb{R}^p} \exp \left( - \frac{1}{2} \bm{x}^\top \big( \tilde{\bm{Q}} + \bm{\Sigma}^{-1} \big) \bm{x} + \big(\bm{\Sigma}^{-1} \bm{\theta} - \tilde{\bm{b}} \big)^\top \bm{x} \right) \mathrm{d}\bm{x} \\
& = \exp \left( - \frac{1}{2} \bm{\theta}^\top \bm{\Sigma}^{-1} \bm{\theta} + \frac{1}{2} \big(\bm{\Sigma}^{-1} \bm{\theta} - \tilde{\bm{b}} \big)^\top \big( \tilde{\bm{Q}} + \bm{\Sigma}^{-1} \big)^{-1} \big(\bm{\Sigma}^{-1} \bm{\theta} - \tilde{\bm{b}} \big) \right) \\
& = \exp \Big( - \frac{1}{2} \bm{\theta}^\top \left( \bm{\Sigma}^{-1} -\bm{\Sigma}^{-1}(\tilde{\bm{Q}} + \bm{\Sigma}^{-1})^{-1}\bm{\Sigma}^{-1} \right) \bm{\theta} 
 - \tilde{\bm{b}}^\top (\tilde{\bm{Q}} + \bm{\Sigma}^{-1})^{-1} \bm{\Sigma}^{-1} \bm{\theta} \Big)\\
& = \exp \Big( - \frac{1}{2} \bm{\theta}^\top \left( \bm{\Sigma}^{-1} -\bm{\Sigma}^{-1}\bm{\Sigma}(\bm{\Sigma}^{-1}-\bm{Q}) \bm{\Sigma} \bm{\Sigma}^{-1} \right) \bm{\theta} 
 - \tilde{\bm{b}}^\top\bm{\Sigma}(\bm{\Sigma}^{-1}-\bm{Q}) \bm{\Sigma}\bm{\Sigma}^{-1} \bm{\theta} \Big)\\
& = \exp \Big( - \frac{1}{2} \bm{\theta}^\top \bm{Q} \bm{\theta} 
 - \bm{b}^\top \bm{\theta} \Big). \qquad (\because \mbox{ definition of }\tilde{\bm{b}})
\end{align*}
Taking the negative logarithm on the both sides leads to \eqref{eq:G-condition}.


\subsection{Proof of Proposition \ref{prop:expectation_condition_for_simple_binom}} 
\label{sec:proof:prop:expectation_condition_for_simple_binom}
\noindent
For $X\sim\mathrm{Bin}(m,\theta)$ and $n\in \mathbb{N}$ such that $n\leq m$, we have
\begin{align}
\mathbb{E}_{p_{x}(\theta)}[X(X-1)\cdots (X-n+1)] = m(m-1) \cdots (m-n+1) \theta^n, \label{eq:binom.identity.1}
\end{align}
because
\begin{align*}
\mathbb{E}_{p_x(\theta)}[X(X-1) \cdots (X-n+1)]
& =  \sum_{k=0}^m k (k-1) \cdots (k-n+1) \binom{m}{k} \theta^k (1-\theta)^{m-k} \\
& = m(m-1) \cdots (m-n+1) \sum_{k=n}^m \binom{m-n}{k-n} \theta^k (1-\theta)^{m-k} \\
& =  m(m-1) \cdots (m-n+1) \theta^n 
\end{align*}
where the last equality holds because $\sum\limits_{k=n}^m \binom{m-n}{k-n} \theta^k (1-\theta)^{m-k}=
\sum\limits_{l=0}^{m-n} \binom{m-n}{l} \theta^{l + n}  (1-\theta)^{m-(l+n)}=
\theta^n \sum\limits_{l=0}^{m-n} \binom{m-n}{l}  \theta^{l}  (1-\theta)^{(m-n) - l}=\theta^n$. 
\eqref{eq:binom.identity.1} implies that 
$\mathbb{E}_{p_{x}(\theta_j)}\big[g_{i,j}(X_j)\big]=\theta_j^{n_{i,j}} \,
(\because \eqref{eq:binom.identity.1} \mbox{ with }n=n_{i,j})$
, and that
\begin{align*}
\mathbb{E}_{p_{\bm{x}}(\bm{\theta})}[G(\bm{X})] 
& = K -  \sum_{i\in[I]}\tilde{a}_{i}\prod_{j\in[p]}\mathbb{E}_{p_{x}(\theta_j)}\big[g_{i,j}(X_j)\big]-\tilde{b} & (\because X_j\mbox{ and }X_h\mbox{ are independent if }j\neq h)\\
& = K -  \sum_{i\in[I]}\tilde{a}_{i}\prod_{j\in[p]}\theta_j^{n_{i,j}}-\tilde{b} & (\because \mathbb{E}_{p_{x}(\theta_j)}\big[g_{i,j}(X_j)\big]=\theta_j^{n_{i,j}})\\
& = K - F(\mathrm{diag}(\bm{u}-\bm{l})\bm{\theta} + \bm{l}).&\Box
\end{align*}

\subsection{Setting of $K$ for \eqref{box_multi_problem2}} \label{sec:Setting.of.K.for.rectangle}
This subsection shows how the upper bound $K$ of $\max\limits_{\bm{x} \in \Xi}G(\bm{x})$ can always be found. 
Since $x_j\in\{ 0,1,\ldots,m_j\}$ and $m_j \geq n_{i,j} > 0$, we have 
$\prod\limits_{k\in[n_{i,j}]}(x_j-k-1)\leq\prod\limits_{k\in[n_{i,j}]}(m_j-k-1)$ 
and thus 
$
 g_{i,j}(x_j) \leq 1
$
by the definition of $g_{i,j}$. This leads to $\tilde{a}_i \prod_{j\in[p]} g_{i,j}(x_j) \leq \tilde{a}_i$ for $\tilde{a}_i \geq 0$ and $\leq 0$ for $\tilde{a}_i  < 0$,  
and 
$
 G(\bm{x}) := 
 \sum\limits_{i\in[I]} \tilde{a}_i \prod\limits_{j\in[p]} g_{i,j}(x_j) + \tilde{b} 
 \leq
 \sum\limits_{i\in[I]} \max \{\tilde{a}_i ,0 \} + \tilde{b},
$ 
so the constant defined as 
$
 K :=
 \sum\limits_{i\in[I]} \max \{\tilde{a}_i ,0 \} + \tilde{b}
$ 
gives the upper bound of the optimal value.


\subsection{Proof of Corollary \ref{cor:update_formula_binomial}}
\label{sec:proof_update_formula_binomial}
The binomial distribution \eqref{eq:binomial.prob} belongs to the exponential family \eqref{eq:expdist} with $\bm{\gamma}(\bm{x})=\bm{x}$, so it is easy to see that for $\bm{\Omega}=\mathrm{diag}(\frac{1}{m_1},\ldots,\frac{1}{m_p})$, we have $\mathbb{E}_{p_{\bm{x}}(\bm{\theta})}[\bm{\Omega} \bm{\gamma}(\bm{X})]=\bm{\theta}$. 
In addition, since $\bm{\eta}(\bm{\theta})=(\ln{\frac{\theta_1}{1-\theta_1}},\ldots,\ln{\frac{\theta_p}{1-\theta_p}})$, we have $ \bm{\Omega} \big( \nabla\bm{\eta}(\bm{\theta}) \big)^{-1}=\mathrm{diag}(\frac{\theta_1(1-\theta_1)}{m_1},\ldots,\frac{\theta_p(1-\theta_p)}{m_p})$, which is equal to $\mathbb{V}_{p_{\bm{x}}(\bm{\theta})}[\bm{\Omega}\bm{\gamma}(\bm{X})]=\mathrm{diag}(\frac{\theta_1(1-\theta_1)}{m_1},\ldots,\frac{\theta_p(1-\theta_p)}{m_p})=\bm{I}_{p_{\bm{x}}}(\bm{\theta})^{-1}$. 
Substituting this into \eqref{eq:natgrad}, we get
$
\bm{\theta}^{(t)} \gets \bm{\theta}^{(t-1)} - \mathrm{diag}\Big(\frac{\theta^{(t-1)}_{1}(1-\theta^{(t-1)}_{1})}{m_1},...,\frac{\theta^{(t-1)}_{p}(1-\theta^{(t-1)}_{p})}{m_p}\Big) \nabla f(\bm{\theta}^{(t-1)}).
$


\noindent
\subsection{Proof of Proposition \ref{prop:EMcond_multinomial}}
\label{sec:proof:prop:expectation_condition_for_simple_multinom}
For $\bm{X}\sim\mathrm{Multi}(m,\bm{\theta})$ and $n_j\in\mathbb{N}$ such that $\sum\limits_{j\in[p]}n_j\leq{m}$, we have
\begin{align}
\mathbb{E}_{p_x(\theta)} \Big[ \prod_{j\in[p]} X_j(X_j-1)\cdots(X_j-n_j+1) \Big]
& =  m(m-1) \cdots \Big(m-\sum_{j\in[p]} n_j+1 \Big) \theta_1^{n_1} \cdots \theta_p^{n_p},
\label{eq:multi.identity.1}
\end{align}
where $X_p=m-\sum\limits_{j\in[p-1]}X_j,\theta_p=1-\sum\limits_{j\in[p-1]}\theta_j$. 
The formula \eqref{eq:multi.identity.1} holds since we have
\begin{align*}
\lefteqn{\mathbb{E}_{p_x(\theta)} \Big[ \prod_{j\in[p]} X_j(X_j-1)\cdots(X_j-n_j+1) \Big]}\\  
& =  \sum_{\bm{x} \in \Xi} \prod_{j\in[p]} x_j(x_j-1)\cdots(x_j-n_j+1) \frac{m!}{x_{1}! \cdots x_{p}!} \theta_{1}^{x_{1}} \cdots \theta_{p}^{x_{p}} \qquad (\Xi = \{ \bm{x} \in \mathbb{N}^{p} ~|~ \bm{x}^\top \bm{1} = m \}) \\
& =  m(m-1) \cdots \Big(m-\sum_{j\in[p]} n_j+1 \Big) \theta_{1}^{n_1} \cdots \theta_{p}^{n_{p}} \sum\limits_{\bm{x} \in \Xi} \frac{(m-\sum_{j\in[p]} n_j)! }{ (x_1-n_1)!\cdots(x_{p}-n_{p})! } \theta_1^{x_1-n_1} \cdots \theta_p^{x_{p}-n_{p}} \\
& \qquad (\Xi = \{ \bm{x} \in \mathbb{N}^{p} ~|~ x_1 \geq n_1,\ldots, x_{p} \geq n_{p}, \bm{x}^\top \bm{1} = m \}) \\
& =  m(m-1) \cdots \Big(m-\sum\limits_{j\in[p]} n_j+1 \Big) \theta_1^{n_1} \cdots \theta_{p}^{n_{p}} 
\underbrace{
\sum_{\bm{x}' \in \Xi'} \frac{(m-\sum_{j\in[p]} n_j)! }{ x_1'!\cdots x_{p}'! }\theta_1^{x_1'} \cdots \theta_{p}^{x_{p}'}
}_{=1} \\
& \qquad (\Xi' = \{ \bm{x}' \in \mathbb{N}^{p} ~|~ {\bm{x}'}^\top \bm{1} = m-\sum_{j\in[p]} n_j \}) \qquad (\because \mbox{change-of-variables for } x_j'=x_j-n_j,\ j\in[p]) \\
& =  m(m-1) \cdots \Big(m-\sum_{j\in[p]} n_j+1 \Big) \theta_1^{n_1} \cdots \theta_{p}^{n_{p}},
\end{align*}
where the last equality holds because the last term of the second line from the bottom is the summation of the probability mass function $p_{\bm{x}'}(\bm{\theta})$ of $\mathrm{Multi}(m-\sum\limits_{j\in[p]}n_j,\bm{\theta})$. \eqref{eq:multi.identity.1} implies that 
$\mathbb{E}_{p_x(\theta)} \big[ g_{i}(\bm{X}) \big]=\prod\limits_{j\in[p]}\theta_j^{n_{i,j}}
\,
(\because \eqref{eq:multi.identity.1} \mbox{ with }n_j=n_{i,j})$,
and that
\begin{align*}
\mathbb{E}_{p_{\bm{x}}(\bm{\theta})}[G(\bm{X})] 
& = K - \sum_{i\in[I]}a_i \mathbb{E}_{p_x(\theta)} \Big[ g_i(\bm{X}) \Big]\\
& = K - \sum_{i\in[I]}a_{i}\prod_{j\in[p]}\theta_j^{n_{i,j}}
 \qquad (\because \mathbb{E}_{p_x(\theta)} \big[ g_i(\bm{X}) \big]=\prod_{j\in[p]}\theta_j^{n_{i,j}})\\
& = K - F(\bm{\theta}). 
\end{align*}


\subsection{Setting of $K$ for \eqref{min_var_problem2}} \label{sec:Setting.of.K.for.simplex}
Similar to Section \ref{sec:Setting.of.K.for.rectangle}, 
this subsection shows how the value $K$ can be given for \eqref{min_var_problem2}. 
Since $x_j\in\{ 0,1,\ldots,m\}$ and $m \geq \sum\limits_{j\in[p]} n_{i,j} > 0$, we have 
$
 \prod\limits_{j\in[p]} \prod\limits_{k\in[n_{i,j}]}(x_j-k+1) \leq \prod\limits_{k\in[\sum_{j\in[p]} n_{i,j}]}(m-k+1)
$ 
and thus 
$
 g_{i}(\bm{x}) \leq 1
$
by the definition of $g_{i}$. This leads to 
$a_i g_{i}(\bm{x}) \leq  a_i$ for $a_i \geq 0$; 
$a_i g_{i}(\bm{x}) \leq 0$ for $a_i  < 0$, 
and 
$
 G(\bm{x}) :=  \sum\limits_{i\in[I]} a_i \prod\limits_{j\in[p]} g_{i,j}(x_j)
 \leq \sum\limits_{i\in[I]}\max \{a_i ,0 \},
$ 
so the constant defined as
$
 K :=  \sum\limits_{i\in[I]} \max \{a_i ,0 \} 
$
gives the upper bound of the optimal value.


\subsection{Proof of Corollary \ref{cor:nat.grad.formula.for.polynom.simplex}} 
\label{sec:proof:cor:nat.grad.formula.for.polynom.simplex}
\noindent 
With 
$\bm{\gamma}(\bm{x})=\bm{x},\bm{\Omega}=\frac{1}{m}\bm{I}\in\mathbb{R}^{(p-1)\times(p-1)}$, 
$\bm{\eta}(\bm{\theta})=(\ln\frac{\theta_1}{1-\bm{1}^\top\bm{\theta}},\ldots,\ln\frac{\theta_{p-1}}{1-\bm{1}^\top\bm{\theta}})^\top$, 
\eqref{eq:unbiased_estimator} is valid. The $(k,l)$-elements of the inverse of the Fisher information matrix is
\begin{align*}
\big(\bm{I}_{p_{\bm{x}}}(\bm{\theta})^{-1}\big)_{k,l} =
 \big(\mathbb{V}_{p_{\bm{x}}(\bm{\theta})}[\bm{\Omega}\bm{\gamma}(\bm{X})]\big)_{k,l} =
\begin{cases}
 \frac{\theta_k(1-\theta_k)}{m}, & k=l\land k,l\in[p-1], \\
 - \frac{\theta_k \theta_l}{m}, & k \neq l\land k,l\in[p-1],
\end{cases}
\end{align*}
which follow from \eqref{fisher_inv} and
\begin{align*}
 \big( \nabla \bm{\eta}(\bm{\theta})^\top \big)_{k,l} = 
\begin{cases}
 \frac{1}{\theta_k} - \frac{1}{1-\bm{1}^\top\bm{\theta}}, & k=l\land k,l\in[p-1], \\
 \frac{1}{1-\bm{1}^\top\bm{\theta}}, & k \neq l\land k,l\in[p-1].
\end{cases}
\end{align*}
we get the update formula \eqref{eq:update.formula.poly.on.unit}.


\subsection{Proof of Theorem \ref{theorem:global_convergence}}\label{sec:proof.Thm:convergence}
The optimality condition for \eqref{eq:EM-step} implies
\begin{align}
\nabla f(\bm{\theta}^{(t)}) + \nabla \bm{\eta}(\bm{\theta}^{(t)})^\top \big(\mathbb{E}_{q_{\bm{x}}(\bm{\theta}^{(t)})}[\bm{\gamma}(\bm{X})]-\mathbb{E}_{q_{\bm{x}}(\bm{\theta}^{(t-1)})}[\bm{\gamma}(\bm{X})]\big) & = \bm{0}.
\label{eq:opt_EM2}
\end{align}
By the definition \eqref{eq:def_eta_ineq} of $\bm{\eta}(\bm{\theta})$, the second term on the left-hand side of \eqref{eq:opt_EM2} is 
\begin{align*}
\sum\limits_{i\in[M]}\sum\limits_{j\in[p]} \frac{r_{j,i}(\mu_j(\bm{\theta}^{(t)})-\mu_j(\bm{\theta}^{(t-1)}))}{g_i(\bm{\theta}^{(t)})} \nabla g_i(\bm{\theta}^{(t)}),
\end{align*}
where we denote $\bm{\mu}(\bm{\theta})=(\mu_1(\bm{\theta}),\ldots,\mu_p(\bm{\theta}))^\top:=\mathbb{E}_{q_{\bm{x}}(\bm{\theta})}[\bm{\gamma}(\bm{X})]$.
With $\nu_i^{(t)}:= \sum_{j\in[p]}r_{j,i}(\mu_j(\bm{\theta}^{(t)})-\mu_j(\bm{\theta}^{(t-1)}))/g_i(\bm{\theta}^{(t)})$, 
we have $\nu_i^{(t)}g_i(\bm{\theta}^{(t)})=\sum_{j\in[p]}r_{j,i}(\mu_j(\bm{\theta}^{(t)})-\mu_j(\bm{\theta}^{(t-1)}))\to 0$ as $t\to \infty$, and 
\eqref{eq:opt_EM2} is represented as
\begin{align} \label{eq:Lagrange_const}
 \nabla
 f(\bm{\theta}^{(t)}) + \sum\limits_{i\in[M]} \nu_i^{(t)} \nabla
 g_i(\bm{\theta}^{(t)}) & = \bm{0}.
\end{align}
If \eqref{eq:ineq_limit} is satisfied, we can let $\nu^*_i:=\lim\limits_{t\to\infty}\nu^{(t)}_i$ and see that 
the KKT condition is fulfilled at the limit $\bm{\theta}=\bm{\theta}^*$:
$
 \nabla f(\bm{\theta}^*) + \sum_{i\in[M]}\nu_i^*\nabla g_i(\bm{\theta}^*) = \bm{0}, ~
 g_i(\bm{\theta}^*) \leq  0, ~\nu_i^* g_i(\bm{\theta}^*) =  0, ~\nu_i^*\geq 0,~i\in[M].
$


\subsection{Fulfillment of the condition \eqref{eq:ineq_limit} of Theorem \ref{theorem:global_convergence}}\label{sec:fulfillment_of_convergence_cond}
\paragraph{\it Case (i): Algorithm \ref{alg:polynom.on.rect}}
The box constraints of \eqref{box_multi_problem2} can be represented as in the form of \eqref{eq:shift-form2_theta} with 
$g_{j}(\bm{\theta})=- \theta_{j},~g_{p+j}(\bm{\theta}) = \theta_{j} - 1,~j\in[p]$ (i.e., $M=2p$). 
It is easy to see that the parameters 
$\bm{\eta}(\bm{\theta})= \big( \ln\frac{\theta_1}{1-\theta_1},\ldots, \ln\frac{\theta_p}{1-\theta_p} \big)^\top$ for the binomial distribution \eqref{eq:binomial.prob} are represented as \eqref{eq:def_eta_ineq} 
with $r_{j,i}:= 1 ~(i=j);~-1~ (i=j+p);~0$ (otherwise) for $i\in[2p],~j\in[p]$. 
The formula \eqref{eq:Expectation_q} then implies that
$\mu_j(\bm{\theta}^{(t-1)})=m_j\theta_j^{(t-1)} - \theta_j^{(t-1)}(1-\theta_j^{(t-1)}) \frac{\partial}{\partial \theta_j}f(\bm{\theta}^{(t-1)}), j\in[p].$ 
To see \eqref{eq:ineq_limit}, observe that for the $i$-th constraint $g_i(\bm{\theta})=-\theta_i,~r_{i,i}=1,~i\in[p]$, we have 
\begin{align*}
\nu^{(t)}_i&:=\frac{r_{i,i}\big(\mu_i(\bm{\theta}^{(t)}_i)-\mu_i(\bm{\theta}^{(t-1)}_i)\big)}{g_i(\bm{\theta}^{(t)})}\\
&=\frac{m_i(\theta^{(t)}_i-\theta^{(t-1)}_i)-\{\theta^{(t)}_i(1-\theta^{(t)}_i)\frac{\partial}{\partial \theta_i}f(\bm{\theta}^{(t)})-\theta^{(t-1)}_i(1-\theta^{(t-1)}_i)\frac{\partial}{\partial \theta_i}f(\bm{\theta}^{(t-1)})\}}{-\theta^{(t)}_i}\\
&=m_i\Big(-1+\frac{\theta^{(t-1)}_i}{\theta^{(t)}_i}\Big)+(1-\theta^{(t)}_i)\frac{\partial}{\partial \theta_i}f(\bm{\theta}^{(t)})-\frac{\theta^{(t-1)}_i}{\theta^{(t)}_i}(1-\theta^{(t-1)}_i)\frac{\partial}{\partial \theta_i}f(\bm{\theta}^{(t-1)})\\
&=(1-\theta^{(t)}_i)\frac{\partial}{\partial \theta_i}f(\bm{\theta}^{(t)}),
\end{align*}
which converges to $\frac{\partial}{\partial \theta_i}f(\bm{\theta}^{*})\geq 0$ as $\theta^{(t)}_i\to 0$ because it would contradict to the monotonicity of $\{f(\bm{\theta}^{(t)})\}$ otherwise. 
The fourth equality follows from the update formula \eqref{eq:update.formula.poly.on.rect}, which implies that 
$m_i\big(\theta^{(t)}_i-\theta^{(t-1)}_i\big)=\theta^{(t)}_i(1-\theta^{(t-1)}_i)\frac{\partial}{\partial \theta_i}f(\bm{\theta}^{(t-1)})$
.
Similarly, for the $(p+i)$-th constraint $g_{i+p}(\bm{\theta})=\theta_i-1,~r_{i+p,i+p}=-1,~i\in[p]$, we have 
$
\nu^{(t)}_i
=-\theta^{(t)}_i\frac{\partial}{\partial \theta_i}f(\bm{\theta}^{(t)}),
$ 
which converges to $-\frac{\partial}{\partial \theta_i}f(\bm{\theta}^{*})\geq 0$ as $\theta^{(t)}_i\to 1$ because it would contradict to the monotonicity of $\{f(\bm{\theta}^{(t)})\}$ otherwise. 
Consequently, in either constraint, we can see that $\lim\limits_{t\to\infty}\nu_i^{(t)}\geq 0$. 
\paragraph{\it Case (ii): Algorithm \ref{alg:polynom.on.simplex}}
The simplex constraint of \eqref{min_var_problem2} corresponds to  
$g_{j}(\bm{\theta})=-\theta_{j},~j\in[p-1]$, and $g_{p}(\bm{\theta})=\bm{1}^\top\bm{\theta}-1$. 
With $r_{j,i}:=1~(i=j);~-1~(i=p);~0$ (otherwise) for $i\in[p],~j\in[p-1]
$, 
the parameters $\bm{\eta}(\bm{\theta})= \big(\ln\frac{\theta_1}{1-\bm{1}^\top\bm{\theta}},\ldots, \ln\frac{\theta_{p-1}}{1-\bm{1}^\top\bm{\theta}} \big)^\top$ for the multinomial distribution are represented in the form of \eqref{eq:def_eta_ineq}. 
The formula \eqref{eq:Expectation_q} then implies that
\begin{align*}
\mu_j(\bm{\theta}^{(t-1)})=m\theta_j^{(t-1)} - \big(\theta_j^{(t-1)}\frac{\partial}{\partial \theta_j}f(\bm{\theta}^{(t-1)}) - \theta_j^{(t-1)}(\bm{\theta}^{(t-1)})^\top\nabla f(\bm{\theta}^{(t-1)})\big),~ j\in[p-1].
\end{align*} 
To see 
\eqref{eq:ineq_limit}, observe that for the $i$-th constraint $g_i(\bm{\theta})=-\theta_i,~r_{i,i}=1,~i\in[p-1]$, we have
\begin{align}
\nu^{(t)}_i
=\frac{\mu_i(\bm{\theta}^{(t)})-\mu_i(\bm{\theta}^{(t-1)})}{-\theta_i^{(t)}} 
=(1-\theta_i^{(t)})\frac{\partial}{\partial \theta_i}f(\bm{\theta}^{(t)}) - \sum_{j \neq i} \theta_j^{(t)}\frac{\partial}{\partial \theta_j}f(\bm{\theta}^{(t)}) \label{eq:nu_unit2}
\end{align}
The second equality follows from the update formula \eqref{eq:update.formula.poly.on.unit}, which implies that 
\begin{align*}
m\big(\theta^{(t)}_i-\theta^{(t-1)}_i\big)=-\theta_i^{(t-1)}\frac{\partial}{\partial \theta_i}f(\bm{\theta}^{(t-1)}) + \theta_i^{(t-1)}(\bm{\theta}^{(t-1)})^\top\nabla f(\bm{\theta}^{(t-1)}). 
\end{align*}
We see from \eqref{eq:nu_unit2} that $\nu^{(t)}_i$ converges to a finite value $\nu^*_i$. To see the sign of $\nu^*$, observe that if $\lim\limits_{t \to \infty}\theta^{(t)}_i = 0$, the first term of \eqref{eq:nu_unit2} is nonnegative because $\frac{\partial}{\partial \theta_i}f(\bm{\theta})$ is continuous and the optimality condition implies $\lim\limits_{t \to \infty} \frac{\partial}{\partial \theta_i}f(\bm{\theta}^{(t)})\geq 0$. We can also show $\lim\limits_{t \to \infty}\theta^{(t)}_i = 1$ similarly. As for the case where $\lim\limits_{t \to \infty}\theta^{(t)}_i \in(0,1)$, it converges to zero since $\lim\limits_{t \to \infty}\frac{\partial}{\partial \theta_i}f(\bm{\theta}^{(t)})=0$ must hold then. 
As for the second term of \eqref{eq:nu_unit2}, we see that $\lim\limits_{t \to \infty}-\theta_j^{(t)}\frac{\partial}{\partial \theta_j}f(\bm{\theta}^{(t)})\geq 0$ for every $j \neq i$, likewise. 
Similarly, for the $p$-th constraint $g_{p}(\bm{\theta})=\bm{1}^\top\bm{\theta}-1$, we have 
\begin{align*}
\nu^{(t)}_p
&
=\frac{-\sum\limits_{j\in[p-1]}\big(\mu_j(\bm{\theta}^{(t)})-\mu_j(\bm{\theta}^{(t-1)})\big)}{\bm{1}^\top\bm{\theta}^{(t)}-1} =- (\bm{\theta}^{(t)})^\top\nabla f(\bm{\theta}^{(t)}) 
\end{align*}
Similar to the case for $i\neq p$, we can show that $\nu^{(t)}_p$ converges to a nonnegative value.

\subsection{Proof of Proposition \ref{prop:G.for.quadratic_problem0}} 
\label{sec:proof:prop:EMcond_for_QP0}
\noindent 
We have
\begin{align*}
\mathbb{E}_{p_{\bm{x}}(\bm{\theta})}[G(\bm{X})]
&= \mathbb{E}_{p_{\bm{x}}(\bm{\theta})}\big[G_1(\bm{X}_1)\,G_2(\bm{X}_2)\big] \\
&= \mathbb{E}_{p_{\bm{x}_1}(\bm{\lambda})}\big[G_1(\bm{X}_1)\big] \, \mathbb{E}_{p_{\bm{x}_2}(\bm{\theta})}\big[G_2(\bm{X}_2)\big] & (\because \bm{X}_1\mbox{ and }\bm{X}_2\mbox{ are independent})\\
&= \exp \Big( - \frac{1}{2} \bm{\theta}^\top \bm{Q} \bm{\theta} - \bm{b}^\top \bm{\theta} \Big),
\end{align*}
where the last equality follows from 
$\mathbb{E}_{p_{\bm{x}_1}(\bm{\lambda})}[G_1(\bm{X}_1)]= \mathbb{E}_{p_{\bm{x}_1}(\bm{\lambda})}[1] = 1,
$ and 
$\mathbb{E}_{p_{\bm{x}_2}(\bm{\theta})}\big[G(\bm{X}_2)\big]
 = \exp ( - \frac{1}{2} \bm{\theta}^\top \bm{Q} \bm{\theta} - \bm{b}^\top \bm{\theta} )
$, which follows from Proof of Proposition \ref{prop:G.for.cvx.q.min}. 
We can also observe the condition $\mathrm{dom}\,p_{\bm{x}}=\{\bm{\theta}\in\mathbb{R}^p\vert\bm{c}-\bm{A}^\top\bm{\theta}\geq\bm{0}\}$. 

\subsection{Proof of Proposition \ref{prop:QPconvergence}}
\label{proof:prop:QPconvergence}
Let $\phi(t)=t\ln t-t+1$ and $H(\bm{x},\bm{y}):=\sum_{j\in[p]}y_j\phi(\frac{x_j}{y_j}) = x_j\ln\frac{x_j}{y_j}+y_j-x_j$, and $L(\bm{\theta}):=(l_1(\bm{\theta}),...,l_p(\bm{\theta}))$ where $l_j(\bm{\theta}):=c_j-\bm{a}_j^\top\bm{\theta},j\in[p]$. 
Here, $H(\bm{x},\bm{y})$ denotes the KL divergence for $\bm{x}\in\mathbb{R}^p_{+}$ relative to $\bm{y}\in\mathbb{R}^p_{++}$.
Observe that Algorithm \ref{alg:EMforQP} is a proximal-like algorithm which iteratively solves 
\begin{align*}
\bm{\theta}^{(t)}\gets\arg\min_{\bm{\theta}}\Big\{
u(\bm{\theta} | \bm{\theta}^{(t-1)})
:=f(\bm{\theta})
+H(L(\bm{\theta}),L(\bm{\theta}^{(t-1)}))
+\frac{1}{2}\|\bm{\theta}-\bm{\theta}^{(t-1)}\|_{\bm{\Sigma}^{-1}-\bm{Q}}^2
\Big\},
\end{align*}
where $\|\bm{\theta}-\bar{\bm{\theta}}\|_{\bm{\Sigma}^{-1}-\bm{Q}}
:=\sqrt{(\bm{\theta}-\bar{\bm{\theta}})^\top(\bm{\Sigma}^{-1}-\bm{Q})(\bm{\theta}-\bar{\bm{\theta}})}$. 
 Following the argument of \cite{teboulle1997convergence} (Theorem 4.3) and combining with a known technique for proving the convergence of proximal point method (with an elliptical norm), we have
\begin{align*}
f(\bm{\theta}^{(t)})-f(\bm{\theta})&\leq H(L(\bm{\theta}),L(\bm{\theta}^{(t-1)}))-H(L(\bm{\theta}),L(\bm{\theta}^{(t)}))+\frac{1}{2}\|\bm{\theta}-\bm{\theta}^{(t-1)}\|_{\bm{\Sigma}^{-1}-\bm{Q}}^2-\frac{1}{2}\|\bm{\theta}-\bm{\theta}^{(t)}\|_{\bm{\Sigma}^{-1}-\bm{Q}}^2
\end{align*}
and consequently, we have for any $\bm{\theta}$
\begin{align*}
f(\bm{\theta}^{(t)})-f(\bm{\theta})&\leq\frac{1}{t}\Big(H(L(\bm{\theta}),L(\bm{\theta}^{(0)}))
+\frac{1}{2}\|\bm{\theta}-\bm{\theta}^{(0)}\|_{\bm{\Sigma}^{-1}-\bm{Q}}^2\Big).
\end{align*}
By substituting $\bm{\theta}=\bm{\theta}^*$, the right-hand side goes to $0$ as $t\to \infty$. 
The second statement can be shown in the same manner as Theorem 4.3(iii) of \cite{teboulle1997convergence} by extending $H(L(\bm{\theta}),L(\bm{\theta}^{(t)}))$ to the composite proximal term, $H(L(\bm{\theta}),L(\bm{\theta}^{(t)}))+\frac{1}{2}\|\bm{\theta}-\bm{\theta}^{(t)}\|_{\bm{\Sigma}^{-1}-\bm{Q}}^2$. 

\section{Quadratic program and linear program}
\label{subsec:QP2}
Consider a quadratic program (QP):
\begin{align*}
({\rm QP.D}) \quad 
\begin{array}{|ll}
\underset{\bm{\theta}_1\in\mathbb{R}^{n},\bm{\theta}_2\in\mathbb{R}^p}{\mbox{maximize}} \quad & \bm{b}^\top\bm{\theta}_1 - \frac{1}{2}\bm{\theta}_2^\top\bm{Q}\bm{\theta}_2 \\
\mbox{subject to} & \bm{A}^\top\bm{\theta}_1-\bm{Q}\bm{\theta}_2\leq\bm{c}, \\
\end{array}
\end{align*}
where $\bm{A}\in\mathbb{R}^{n \times p},\bm{b} \in \mathbb{R}^n,\bm{c}\in\mathbb{R}^p$, and $\bm{Q}\in\mathbb{S}^p$. 
(QP.D) includes a linear program as a special case by setting $\bm{Q}=\bm{O}$ and ignoring $\bm{\theta}_2$. 
(QP.D) is the Lagrangian dual to another convex QP:
\begin{align*}
({\rm QP.P}) \quad 
\begin{array}{|ll}
\underset{\bm{\theta}_2\in\mathbb{R}^p}{\mbox{minimize}} \quad & \frac{1}{2}\bm{\theta}_2^\top\bm{Q}\bm{\theta}_2+\bm{c}^\top\bm{\theta}_2 \\
\mbox{subject to} & \bm{A}\bm{\theta}_2=\bm{b},\phantom{\bm{\theta}_2\geq\bm{0}}\leftarrow\bm{\theta}_1\\
                         & \bm{\theta}_2\geq\bm{0},\phantom{\bm{A}\bm{\theta}_2=\bm{b}}\leftarrow\bm{\lambda}\geq\bm{0},
\end{array}
\end{align*}
where the correspondence between constraints and Lagrange multipliers, $\bm{\theta}_1$ and $\bm{\lambda}$, is indicated.  
We assume that $n\leq p$ and the rank of $\bm{A}$ is $n$, so that there exists a vector $\bm{\xi}$ satisfying $\bm{A}\bm{\xi}=\bm{b}$. 
For $S>\max\big\{\max\limits_{j\in [p]}\{\xi_j\},0\big\}$, define $\hat{\bm{\xi}}:=\frac{1}{S}\bm{\xi},\hat{\bm{Q}}:=\frac{1}{S}\bm{Q}$, and 
another convex QP:
\begin{align*}
({\rm QP.D}') \quad 
\begin{array}{|ll}
\underset{\bm{\theta}_1\in\mathbb{R}^{n},\bm{\theta}_2\in\mathbb{R}^p}{\mbox{minimize}} \quad & f(\bm{\theta}_1,\bm{\theta}_2)=
\hat{\bm{\xi}}^\top(\bm{c}-\bm{A}^\top\bm{\theta}_1)+\frac{1}{2} \bm{\theta}_2^\top\hat{\bm{Q}}\bm{\theta}_2 \\
\mbox{subject to} & \bm{c}-\bm{A}^\top\bm{\theta}_1+\bm{Q}\bm{\theta}_2\geq\bm{0}.
\end{array}
\end{align*}
We see that ({\rm QP.D}$'$) is equivalent to ({\rm QP.D}) since $f(\bm{\theta}_1,\bm{\theta}_2)= 
\hat{\bm{\xi}}^\top(\bm{c}-\bm{A}^\top\bm{\theta}_1)+\frac{1}{2} \bm{\theta}_2^\top\hat{\bm{Q}}\bm{\theta}_2=-\frac{1}{S}\left( \bm{b}^\top \bm{\theta}_1 - \frac{1}{2} \bm{\theta}_2^\top \bm{Q} \bm{\theta}_2 \right) + \frac{1}{S}\bm{\xi}^\top\bm{c}=-\frac{1}{S}``$objective of ({\rm QP.D})
''$+$``constant.'' 

To apply EM algorithm to (QP.D$'$), consider $2p$ random variables $\bm{X}:=(\bm{X}_1,\bm{X}_2)$, where 
$\bm{X}_1$ denotes $p$ independent Poisson random variables, $\bm{X}_2$ follows $p$-dimensional normal distribution $\bm{\mathrm{N}}(\bm{\theta}_2,\bm{\Sigma})$ for $\bm{\Sigma}\in\bm{S}^p_{++}$, and every pair of $X_{1,i}$ and $X_{2,j}$ is independent for any $i,j\in [p]$. 
The likelihood-like function for $\bm{X}=(\bm{x}_1,\bm{x}_2)$ is proportional to
\begin{align}
p_{\bm{x}}(\bm{\theta}_1,\bm{\theta}_2) & = p_{\bm{x}_1}(\bm{\lambda})\,p_{\bm{x}_2}(\bm{\theta}_2)
= p_{\bm{x}_1}(\bm{c}-\bm{A}^\top\bm{\theta}_1+\bm{Q} \bm{\theta}_2)\,p_{\bm{x}_2}(\bm{\theta}_2)
= \hat{p}_{\bm{x}_1}(\bm{\theta}_1,\bm{\theta}_2)\,p_{\bm{x}_2}(\bm{\theta}_2)
,
\label{eq:Poisson+Normal}
\end{align}
where 
\begin{align*}
p_{\bm{x}_1}(\bm{\lambda}) 
&:=\prod_{j\in[p]}\frac{\lambda_{j}^{x_{1,j}}\exp(-\lambda_{j})}{x_{1,j}!} 
\\
&=\hat{p}_{\bm{x}_1}(\bm{\theta}_1,\bm{\theta}_2) 
:=\prod_{j\in[p]}\frac{\big(c_j - \bm{a}_j^\top \bm{\theta}_1 + \bm{q}_j^\top \bm{\theta}_2 \big)^{x_{1j}}\exp\big(-(c_j - \bm{a}_j^\top \bm{\theta}_1 + \bm{q}_j^\top \bm{\theta}_2)\big)}{x_{1,j}!}, \\
p_{\bm{x}_2}(\bm{\theta}_2) 
& :=\frac{1}{(2\pi)^{\frac{p}{2}} \sqrt{|\bm{\Sigma}|}} \exp\Big(- \frac{1}{2}(\bm{x}_2-\bm{\theta}_2)^\top\bm{\Sigma}^{-1} (\bm{x}_2-\bm{\theta}_2) \Big).
\end{align*}
It is easy to see that \eqref{eq:Poisson+Normal} is a member of the exponential family.
\begin{proposition}
\label{prop:G.for.quadratic_problem}
Suppose that $\bm{\Sigma}^{-1}\succ\bm{Q}$. If $G$ is defined by 
$
G(\bm{x}_1,\bm{x}_2) = G_1(\bm{x}_1)G_2(\bm{x}_2),
$ 
where
\begin{align*}
G_1(\bm{x}_1) &:= \exp \Big( \sum_{j\in[p]} x_{1,j} \ln(1-\hat{\xi}_j) \Big), \quad
G_2(\bm{x}_2) := \sqrt{\frac{|\bm{\Sigma}|}{|\tilde{\bm{Q}} + \bm{\Sigma}^{-1}|}} \exp \Big( - \frac{1}{2}\bm{x}_2^\top \tilde{\bm{Q}}\bm{x}_2 - \tilde{\bm{b}}^\top \bm{x}_2 \Big),\\
 \tilde{\bm{Q}} &:= \bm{\Sigma}^{-1}(\bm{\Sigma}^{-1}-\hat{\bm{Q}})^{-1} \bm{\Sigma}^{-1} - \bm{\Sigma}^{-1}, \quad\mbox{and}\quad
 \tilde{\bm{b}} := -\bm{\Sigma}^{-1}(\bm{\Sigma}^{-1}-\hat{\bm{Q}})^{-1} \hat{\bm{Q}}\hat{\bm{\xi}},
\end{align*}
the condition \eqref{eq:G-condition_ext} for (QP.D$'$) is fulfilled 
with $
p_{\bm{x}}(\bm{\theta}_1,\bm{\theta}_2)$ defined by \eqref{eq:Poisson+Normal}, 
$f(\bm{\theta}_1,\bm{\theta}_2) = \hat{\bm{\xi}}^\top(\bm{c}-\bm{A}^\top\bm{\theta}_1
) 
+\frac{1}{2} \bm{\theta}_2^\top\hat{\bm{Q}}\bm{\theta}_2$, and $\Theta=\mathbb{R}^p_+$.
\end{proposition}
%
\noindent
\paragraph{Proof.}
Recalling the moment generating function (MGF) of Poisson distribution: 
$\mathbb{E}_{p_{x_{1,j}}(\lambda_j)}[\exp(t X_{1,j})] = \exp \big( \lambda_j (\exp(t)-1) \big)$,  
we have
\begin{align}
\mathbb{E}_{p_{\bm{x}_1}(\bm{\lambda})}[G_1(\bm{X}_1)] 
& = \mathbb{E}_{p_{\bm{x}_1}(\bm{\lambda})}\Big[\exp \Big( \sum_{j\in[p]} X_{1,j} \ln(1-\hat{\xi}_j) \Big)\Big] \nonumber\\
& = \prod_{j\in[p]} \mathbb{E}_{p_{x_{1,j}}(\lambda_j)}\big[\exp \big( X_{1,j} \ln(1-\hat{\xi}_j) \big)\big] \nonumber& (\because X_{1,j}
\mbox{ are independent}) \nonumber\\
& = \prod_{j\in[p]} \exp \big( - \hat{\xi}_j \lambda_j \big) & (\because \mbox{ MGF with }t=\ln(1-\hat{\xi}_j))\nonumber\\
& = \exp \big( - \hat{\bm{\xi}}^\top \bm{\lambda} \big).\label{eq:poisson_leads_to_linear_logexpect}
\end{align}
Putting this and Propositions \ref{prop:G.for.cvx.q.min} together, we have
\begin{align*}
\mathbb{E}_{p_{\bm{x}}(\bm{\theta}_1,\bm{\theta}_2)}[G(\bm{X})]
&= \mathbb{E}_{p_{\bm{x}}(\bm{\theta}_1,\bm{\theta}_2)}\big[G_1(\bm{X}_1)\,G_2(\bm{X}_2)\big] \\
&= \mathbb{E}_{p_{\bm{x}_1}(\bm{\lambda})}\big[G_1(\bm{X}_1)\big] \, \mathbb{E}_{p_{\bm{x}_2}(\bm{\theta}_2)}\big[G_2(\bm{X}_2)\big] & (\because \bm{X}_1\mbox{ and }\bm{X}_2\mbox{ are independent})\\
&= \exp \big( - \hat{\bm{\xi}}^\top \bm{\lambda} \big) \exp \Big( - \big( -\hat{\bm{\xi}}^\top\bm{Q} \bm{\theta}_2 + \frac{1}{2} \bm{\theta}_2^\top \hat{\bm{Q}} \bm{\theta}_2 \big) \Big) & (\because \mbox{ Proposition }\ref{prop:G.for.cvx.q.min}\mbox{ and }\eqref{eq:poisson_leads_to_linear_logexpect})
\\
&= \exp \Big( - \big( \hat{\bm{\xi}}^\top \bm{\lambda}  -\hat{\bm{\xi}}^\top\bm{Q} \bm{\theta}_2 + \frac{1}{2} \bm{\theta}_2^\top \hat{\bm{Q}} \bm{\theta}_2 \big) \Big),
\end{align*}
where $\bm{\lambda}=\bm{c}-\bm{A}^\top\bm{\theta}_1+\bm{Q} \bm{\theta}_2$.\hfill$\Box$


For $\bm{\theta},\bar{\bm{\theta}}\in
\mathrm{int}\{\bm{\theta}=(\bm{\theta}_1,\bm{\theta}_2)\vert \bm{A}^\top\bm{\theta}_1-\bm{Q} \bm{\theta}_2\leq\bm{c}\}$, 
the surrogate 
$u(\bm{\theta} | \bar{\bm{\theta}})$ is 
a convex function of the form:
\begin{align*}
u(\bm{\theta} | \bar{\bm{\theta}}) & = -\sum\limits_{j\in[p]} \Big\{(1-\hat{\xi}_j)(c_j-\bm{a}_j^\top\bar{\bm{\theta}}_1+\bm{q}_j^\top\bar{\bm{\theta}}_2) 
\ln\frac{c_j-\bm{a}_j^\top\bm{\theta}_1+\bm{q}_j^\top\bm{\theta}_2}{c_j-\bm{a}_j^\top\bar{\bm{\theta}}_1+\bm{q}_j^\top\bar{\bm{\theta}}_2}+\bm{a}_j^\top(\bm{\theta}_1-\bar{\bm{\theta}}_1) -\bm{q}_j^\top(\bm{\theta}_2-\bar{\bm{\theta}}_2)\Big\} \\
& \qquad - \big\{ 
(\bm{\Sigma}^{-1}-\bm{Q})\bar{\bm{\theta}}_2-\bm{Q}\hat{\bm{\xi}}
\big\}^\top
\big(\bm{\theta}_2-\bar{\bm{\theta}}_2 \big) + \frac{1}{2}\bm{\theta}_2^\top\bm{\Sigma}^{-1}\bm{\theta}_2 - \frac{1}{2}\bar{\bm{\theta}}_2^\top\bm{\Sigma}^{-1}\bar{\bm{\theta}}_2 + f(\bar{\bm{\theta}}),
\end{align*}
and the EM algorithm for ({\rm QP.D}) is described in Algorithm \ref{alg:EMforQP}.
\begin{algorithm}[tb]
\caption{EM algorithm for quadratic program (QP.D$'$)}
\label{alg:EMforQP}
\begin{algorithmic}[1]
\State 
Let $\bm{\xi}$ such that $\bm{A}\bm{\xi}=\bm{b}$, 
 $S>\max\big\{\max\limits_{j\in [p]}\{\xi_j\},0\big\}$,
and $(\bm{\theta}_1^{(0)},\bm{\theta}_2^{(0)})\in\mathbb{R}^n\times\mathbb{R}^{p}$, 
and set $t\gets 1$. 
\Repeat
\State
For 
$\mu^{(t-1)}_j:= (1-\hat{\xi}_j)(c_j-\bm{a}_j^\top\bm{\theta}_1^{(t-1)}+\bm{q}_j^\top\bm{\theta}_2^{(t-1)})$ and 
$\bm{\nu}^{(t-1)}:= 
\bm{Q}\hat{\bm{\xi}}-(\bm{\Sigma}^{-1}-\bm{Q})\bm{\theta}^{(t-1)}_2
$, 
solve
\begin{align}
(\bm{\theta}_1^{(t)},\bm{\theta}_2^{(t)})\gets
\arg\min_{(\bm{\theta}_1,\bm{\theta}_2)}
\left\{
\begin{array}{r}
-\sum\limits_{j\in[p]} \Big( \bm{a}_j^\top\bm{\theta}_1 -\bm{q}_j^\top\bm{\theta}_2  + \mu^{(t-1)}_j \ln(c_j-\bm{a}_j^\top\bm{\theta}_1+\bm{q}_j^\top\bm{\theta}_2 \Big)\\
+ \frac{1}{2}\bm{\theta}_2^\top\bm{\Sigma}^{-1}\bm{\theta}_2+(\bm{\nu}^{(t-1)})^\top\bm{\theta}_2
\end{array}
\right\},
\label{eq:M-step_for_QP_2}
\end{align}
\State and let $t\gets t+1$.
\Until{a termination condition is fulfilled}
\end{algorithmic}
\end{algorithm}

The convergence of Algorithm \ref{alg:EMforQP} can be shown based on \cite{teboulle1997convergence} in a similar way to Proposition \ref{prop:QPconvergence} by setting $l_j(\bm{\theta}):=(1-\hat{\xi}_j)(c_j-\bm{a}_j^\top\bm{\theta}_1+\bm{q}_j^\top\bm{\theta}_2),j\in[p]$ and replacing the proximal like algorithm with
$
\bm{\theta}^{(t)}\gets\arg\min_{\bm{\theta}}\Big\{
u(\bm{\theta} | \bm{\theta}^{(t-1)})
:=f(\bm{\theta})
+H(L(\bm{\theta}),L(\bm{\theta}^{(t-1)}))
+\frac{1}{2}\|\bm{\theta}_2-\bm{\theta}_2^{(t-1)}\|_{\bm{\Sigma}^{-1}-\bm{Q}}^2
\Big\},
$ 
where $\|\bm{\theta}_2-\bar{\bm{\theta}}_2\|_{\bm{\Sigma}^{-1}-\bm{Q}}
:=\sqrt{(\bm{\theta}_2-\bar{\bm{\theta}}_2)^\top(\bm{\Sigma}^{-1}-\bm{Q})(\bm{\theta}_2-\bar{\bm{\theta}}_2)}$.

It is interesting to see that $\bm{\theta}_2$ is multiplied by $\bm{Q}$ in (QP.D) or (QP.D$'$), and the 
normal distribution is not needed when we consider LP via (QP.D) or (QP.D$'$). 
We can modify Proposition \ref{prop:G.for.quadratic_problem} to use only $p$ Poisson distributions to approach such LP:
\begin{align*}
({\rm LP.D}') \quad 
\begin{array}{|ll}
\underset{\bm{\theta}\in \mathbb{R}^n}{\mbox{minimize}} \quad & f(\bm{\theta}):= 
\hat{\bm{\xi}}^\top(\bm{c}-\bm{A}^\top\bm{\theta}) \\
\mbox{subject to} & \bm{c}-\bm{A}^\top\bm{\theta}\geq\bm{0}.
\end{array}
\end{align*}

\begin{corollary}
\label{cor:G.for.LP}
For $p_{\bm{x}}(\bm{\lambda}) 
=\hat{p}_{\bm{x}}(\bm{\theta}) 
:=\prod_{j\in[p]}\frac{\big(c_j - \bm{a}_j^\top \bm{\theta} \big)^{x_{j}}\exp\big(-(c_j - \bm{a}_j^\top \bm{\theta})\big)}{x_{j}!}$ and $G$ is defined by 
$G(\bm{x}) := \exp \Big( \sum_{j\in[p]} x_{j} \ln(1-\hat{\xi}_j) \Big)$, 
the condition \eqref{eq:G-condition_ext} for (LP.D$'$) is fulfilled with $
p_{\bm{x}}(\bm{c}-\bm{A}^\top\bm{\theta})$ and $f(\bm{\theta}) = \hat{\bm{\xi}}^\top(\bm{c}-\bm{A}^\top\bm{\theta})$ and $\Theta=\mathbb{R}^p_+$.
\end{corollary}


\section{EM gradient algorithms for the extended formulations}\label{sec:GEM}
This section collectively states EM gradient algorithms \cite{lange1995gradient} for the problems presented in Sec.~\ref{sec:further_examples}. 

Algorithm \ref{alg:NewtonGEM} shows an EM gradient algorithm for the subproblem \eqref{eq:EM-step}. 
\begin{algorithm}[tb]
\caption{An EM gradient algorithm for \eqref{eq:shift-form2}}
\label{alg:NewtonGEM}
\begin{algorithmic}[1]
\State Let $\bm{\theta}^{(0)}$ be an interior point of the feasible region  of \eqref{eq:shift-form2}, 
$\bar{f}\gets f(\bm{\theta}^{(0)}) (= u(\bm{\theta}^{(0)}\vert\bm{\theta}^{(0)}))$ and $\beta\in(0,1)$.  
\Repeat
\State Set $\alpha\gets 1$, and compute
\begin{align}
\hat{\bm{\theta}} \gets \bm{\theta}^{(t-1)} - \alpha \left(\nabla^2 u(\bm{\theta} | \bm{\theta}^{(t-1)})\big\vert_{\bm{\theta}=\bm{\theta}^{(t-1)}}\right)^{-1}\left.\nabla u(\bm{\theta} | \bm{\theta}^{(t-1)})\right\vert_{\bm{\theta}=\bm{\theta}^{(t-1)}}
.
\label{eq:NewtonGEM}
\end{align}
\While{$\hat{\bm{\theta}}\not\in\mathrm{int}\Theta$ or $f(\hat{\bm{\theta}})>\bar{f}$}
\State $\alpha\gets\alpha\beta$ and update $\hat{\bm{\theta}}$ via \eqref{eq:NewtonGEM}, and set 
$\bar{f}\gets f(\hat{\bm{\theta}})$ if $\hat{\bm{\theta}}\in\mathrm{int}\Theta$ and $f(\hat{\bm{\theta}})<\bar{f}$.
\EndWhile
\State Set $\bm{\theta}^{(t)}\gets\hat{\bm{\theta}}$ 
 and $t\gets t+1$.
\Until{a termination condition is satisfied.}
\end{algorithmic}
\end{algorithm}
In the vanilla EM gradient, the Newton update can go out of the domain of $u$, so 
a backtracking procedure is employed for the inner loop (lines 4 to 6) to ensure the update to be strictly feasible to the original problem and the objective value to decrease. 
Note that due to the differentiability of $u$ and $f$ and the first-order surrogate condition, this inner loop will end in a finite number of backtracking. 

\end{appendices}

\end{document}